\documentclass[a4paper,11pt]{article}
\usepackage[latin1]{inputenc}

\pdfoutput=1

\usepackage{amsmath,amsfonts,amssymb} 
\usepackage{graphicx}
\usepackage{enumerate}

\usepackage{tikz}
\usetikzlibrary{matrix}

\usepackage{special2}

\newcommand{\theofont}{\sf \large }
\newtheorem{Def}{\theofont Definition}
\newtheorem{Rem}{\theofont Remark}[subsection]
\newtheorem{Example}[Rem]{\theofont Example}
\newtheorem{Lem}[Rem]{\theofont Lemma}
\newtheorem{Prop}[Rem]{\theofont Proposition}

\newtheorem{Theo}{\theofont Theorem}

\newtheorem{Cor}[Rem]{\theofont Corollary}
\newcommand{\PR}{{\vspace*{-0.2cm}\theofont{\em{Proof.}} } }
\newcommand{\PRnotheo}{{\vspace*{0.2cm}\noindent\theofont{\em{Proof.}} } }
\newcommand{\EPR}{\hfill $\Box$ \vspace*{0.3cm}}
\newcommand{\EEX}{\hfill $\Diamond$ \end{Example} \vspace*{0.3cm}}
\newcommand{\EQ}{\begin{equation}}
\newcommand{\EEQ}{\end{equation}}
\newcommand{\smbf} {\bf \small}

\newcommand{\Lang}{L}
\newcommand{\SP}{SP_{\Lang}}

\newcommand{\Tpreiet}{f_\Lang} 
\newcommand{\Ipreiet}{X_\Lang} 
\newcommand{\Tiet}{T_\Lang} 
\newcommand{\Tiettm}{T_{\Lang_{tm}}} 
\newcommand{\TietP}{T_{\Lang_{P}}}
\newcommand{\TietPA}{T_{\Lang_{P_A}}}
 
\newcommand{\Ik}{I^{(k)}} 
\newcommand{\Jk}{J^{(k)}} 
\newcommand{\Ikiet}{I^{(k)}_\Lang} 
\newcommand{\Iaiet}{I_{\Lang,a}}
\newcommand{\Bacc}{B_{acc}} 
 
\newcommand{\Disc}{{\cal D}}
\newcommand{\Discacc}{{\cal D}_{acc}}
\newcommand{\Discaccr}{{\cal D}_{acc,r}}

\newcommand{\Disciet}{{\cal D}_\Lang}
\newcommand{\Discietacc}{{\cal D}_{\Lang,acc}}
\newcommand{\Discietaccr}{{\cal D}_{\Lang,acc,r}}

\newcommand{\Biet}{B_\Lang}  
\newcommand{\Bietacc}{B_{\Lang,acc}} 
\newcommand{\Bietaccr}{B_{\Lang,acc,r}}

\newcommand{\Bietm}{B_{\Lang_,m}}

\newcommand{\Iietmi}{I_{\Lang,m,i}} 

\newcommand{\Baccone}{B_{acc}^1} 
\newcommand{\Bacctwo}{B_{acc}^2} 
\newcommand{\Bietaccone}{B_{\Lang,acc}^1} 
\newcommand{\Bietacctwo}{B_{\Lang,acc}^2} 

\newcommand{\vmaxpref}{v_{\mbox{\it \scriptsize pref}}}
\newcommand{\vmaxprefl}{v_{\mbox{\it \scriptsize pref},l}}
\newcommand{\vmaxprefr}{v_{\mbox{\it \scriptsize pref},r}}
\newcommand{\vmaxprefgen}{v_{\mbox{\it \scriptsize pref}}}

\newcommand{\Tseq}{{\cal T}}

\begin{document}

\author{\small
  Luis-Miguel Lopez\\
  {\footnotesize Tokyo University of Social Welfare}\\ 
  \and 
  \small Philippe Narbel\\
  {\footnotesize LaBRI, University of Bordeaux}\\
}

\title{\vspace*{-0.4cm}
       \Large \sc Infinite Interval Exchange Transformations \\ from Shifts}
\date{\vspace*{0.5cm}
}

\maketitle

\vspace*{0.2cm}
\noindent
{\small {\sc Abstract.} {\em We show that minimal shifts with zero topological
    entropy are topologically conjugate to interval exchange transformations,
    generally infinite.  When these shifts have linear factor complexity (linear
    block growth), the conjugate interval exchanges are proved to satisfy strong
    finiteness properties.
    }
}  \\

\noindent
{\footnotesize {\sc Keywords:} interval exchange transformation; symbolic
  dynamics; shift; factor complexity; topological entropy; infinite special
  word.}
\vspace*{0.8cm}
\section{Introduction}

\vbox{
Interval exchange transformations (IETs) are maps over $I=[0,1)$ which can be
  defined as permutations of intervals partitioning $I$.  In the finite case,
  these maps are just piecewise isometries of $I$ onto itself, or equivalently,
  injective maps preserving measure and having only a finite number of
  discontinuities.
They happen to be a fundamental notion in dynamical systems and ergodic
theory~\cite{CFS82,Man87,HK02}.
Another important notion in the same context
is the {shift}~$\sigma$ on the set $A^\mathbb{N}$ of infinite words
(symbolic sequences) over a finite alphabet~$A$:
this simple continuous map consists in erasing the first letter of its argument.
The pair $({A}^\mathbb{N},\sigma)$ forms a basic topological dynamical
system, where~${A}^\mathbb{N}$ is then called the {\em full shift} over~$A$;  
if $L$ is a closed $\sigma$-invariant subset of ${A}^\mathbb{N}$, the pair
$(L,\sigma)$ also induces such a system,
where $L$ is then just called a {\em shift}~\cite{LM95,Kit98}.
The topological entropy of a shift $L$ depends on the factor complexity
(block growth) of~$\Lang$~\cite{MH38,Par66,CN10}, 
i.e.~the map $p_\Lang$ on~$\mathbb{N}^*$ giving for each~$n$ the number of
distinct length-$n$ factors (subblocks) occurring in $L$'s infinite words.
Now, a known general relationship between all the above concepts is the
following: the support~$I$ of a finite IET can be embedded as a subset into a
measured compact space in such a way that the IET extends~to
a measure-preserving continuous map, 
whose natural symbolic conjugate is a shift~$L$ with a factor complexity bounded
by an affine function (thus $L$ has zero topological entropy)~\cite{Kea75}.
The main idea of this paper is to study this relationship the other way
around, that is, starting from shifts to build topologically conjugate IETs,  
generally with infinitely many discontinuities, 
and using factor complexity to determine families of these IETs. 
}

Thus, given a shift $\Lang$ we first describe an ordered measured~compact
space~$\Ipreiet$ containing $I$, and a continuous self-map $\Tpreiet$ on
$\Ipreiet$
in such a way that $(\Lang, \sigma)$ is topologically conjugate to $(\Ipreiet,
\Tpreiet)$.
Next, we show how~$\Tpreiet$ can be seen as an extension over $\Ipreiet$ of an
IET $\Tiet$ over~$I$.  We shall mostly consider {aperiodic minimal shifts},
i.e.~shifts with no periodic word and containing no proper shift.
With this respect, the first main result we prove is~the~following: 

\begin{Theo}\label{main1}
Let $\Lang$ be a measured minimal aperiodic shift with zero topological
entropy.
Then $(\Lang, \sigma)$ is topologically~conjugate to $(\Ipreiet,\Tpreiet)$,
where~$\Tpreiet$ is the extension over~$\Ipreiet$ of an IET $\Tiet$ on~$I
\subset \Ipreiet$.
\end{Theo}

\noindent
This result is to be put into perspective with the fact that aperiodic measure
preserving transformations of a Lebesgue space are known to be isomorphic to 
infinite IETs~\cite{AOW85}:
Besides the fact it defines a tight relationship between shifts and
transformations over $[0,1)$, Theorem~\ref{main1} is about topological
  isomorphisms --~i.e.~homeomorphisms~--, and not only about isomorphisms
  --~i.e.~measure-preserving maps.
Another comparison point is that the conjugate IETs given by Theorem~\ref{main1}
can have up to a null measure infinite set of discontinuities, while being
always piecewise increasing.
Also, the construction behind these conjugate maps has the following consistency
property: starting with a piecewise increasing~IET~$T$, and coding its dynamics
after a finite monotonicity partition of~$I$, a shift is obtained to which
Theorem~\ref{main1} applies, yielding an IET coinciding with~$T$
(see Proposition~\ref{iet_gives_the_same_iet_prop} further).

The paper focuses next on aperiodic minimal shifts $\Lang$ with linear
complexity~\cite{Fog02,CN10}, i.e.~such that $p_L(n)=O(n)$, forming a family of
shifts with zero~topological entropy, which contains all the simplest
non-trivial ones.
In particular, these linear complexity shifts include the conjugate shifts
obtained from the natural coding of finite IETs,
but also e.g.~the shifts associated with primitive substitutions~\cite{Que10}.
Now, the simplicity behind linear complexity reflects in the conjugate
IETs.  As a matter of fact, in this case, a conjugate IET given by
Theorem~\ref{main1} has three properties: in addition to be piecewise
increasing, its discontinuities may accumulate only at a finite set, and all
these discontinuities belong to only finitely many distinct iterates (full
orbits) of the IET.
We call {\em almost finite} an~IET satisfying these properties, and the second
main result we prove is the following variation of Theorem~\ref{main1}:

\begin{Theo}\label{main2}
Let $\Lang$ be a measured minimal aperiodic shift with linear complexity. 
Then $(\Lang, \sigma)$ is topologically conjugate to $(\Ipreiet,\Tpreiet)$,
where $\Tpreiet$ is the extension over~$\Ipreiet$ of an almost finite IET
$\Tiet$ on $I \subset \Ipreiet$.
\end{Theo}

The above theorems can be constructive. In particular, we give a technique to
exhibit $\sigma$-invariant measures on the shifts $\Lang$ for which the
conjugate IETs can be explicitly built
(see~Proposition~\ref{approx_prop}, and the examples illustrating it).

\section{Basic Definitions}
\label{basics_defs}

\subsection{Interval Exchange Transformations}
\label{iet_subsec}

\begin{Def}\label{iet_def}
An orientation preserving {\smbf finite interval exchange
  transformation}\footnote{In this paper, it is sufficient to consider the
  orientation preserving case only, i.e.~interval exchange transformations
  with~$(+1)$~slopes, and no $(-1)$~slopes as in the general
  case~\cite{Man87}.}  (IET)
is a map $T: I\rightarrow I$, where $I=[0,1)$, with: 
\begin{itemize}\itemsep=0cm
\vspace*{-0.1cm}
\item A finite set $\{x_i\}_{i=0,...,m-1}$ of discontinuities, denoted by
  $\Disc$, with $x_0=0 < x_1 <...< x_{m-1}< x_m = 1$, which determines an
  ordered partition $\bigsqcup_{i=0}^{m-1} I_i$ of~$I$ formed by the right-open
  intervals $I_i=[x_{i}, x_{i+1})$,
\end{itemize}
\vspace*{-0.1cm}
such that: 
\vspace*{-0.1cm}
\begin{itemize}
\item $T$ is injective, and $T$ is a translation on each interval $I_i$,
  i.e.~for each $i=0,...,m-1$, there exists $k_i \in \mathbb{R}$ such that for
  all $x \in I_i$, $T(x) = x + k_i$.
\end{itemize}
\vspace*{-0.2cm}
\end{Def}

\noindent 
Such a finite IET $T$ is just a piecewise order-preserving isometry of $I$.  It
is right-continuous, and measure preserving, being injective and with
derivative~1 on~$I\setminus \Disc$.
A usual way of interpreting~$T$ --~giving its name to the notion~-- is the
following: the intervals $I_i$, as ordered components of the partition
$\bigsqcup_{i=0}^{m-1} I_i$ of $I$,
are permuted 
in the image of~$T$ so as to form another partition $\bigsqcup_{i=0}^{m-1}
I_{\pi(i)}$ of~$I$, where $\pi$ is the permutation over $\{0,...,m-1\}$ induced
by the $k_i$'s determining $T$.
An IET can be indeed also defined by a pair~$(\lambda,\pi)$, where~$\lambda$ is
the
vector of the $I_i$'s lengths and $\pi$ is a permutation of the~$I_i$'s.
Note that given an IET, the partition of $I$ can be refined by using any finite
set $B \supset \Disc$ of points, and such a refinement can be used as well in
Definition~\ref{iet_def} by replacing $\Disc$ by $B$.  The points determining a
specific partition of an IET~$T$ are then usually called the {\smbf separation
  points} of~$T$.

The above definition of an IET can be extended to include infinity as follows:

\begin{Def}\label{iiet_def}
An orientation preserving {\smbf infinite interval exchange transformation} is a
map $T: I\rightarrow I$, where $I=[0,1)$, with:
\vspace*{-0.1cm}
\begin{itemize}\itemsep=0cm
  \label{infinite-iet-def}
  \item An infinite set $Y$ of discontinuities such that the closure of $Y$ in
    $I$,
    denoted by~$\Disc$, has null measure and determines an ordered partition
    of~$I$ formed by:\\
    (i)~The right-open intervals $[x, x')$ with $x \in \Disc$, $x' \in \Disc \cup
      \{1\}$, $[x, x') \cap Y \subseteq \{x\}$;
    (ii)~The set $\Discaccr$ of the accumulation points of $\Disc$ from the
        right,
\end{itemize}
\vspace*{-0.1cm}
such that: 
\vspace*{-0.1cm}
\begin{itemize}
  \item $T$ is injective on $I\setminus \Discaccr$, right-continuous on $I$, and
    $T$ is a translation on each right-open interval of the above
    partition\footnote{Infinite IETs are considered in different ways in the
      literature depending on the authors, e.g.~in~\cite{AOW85}, the
      discontinuities accumulate only at $1$; in~\cite{Hoo15}, the
      discontinuities are countable and the partitioned interval may have
      infinite length.}.
\end{itemize}
\vspace*{-0.1cm}
\end{Def}

\noindent
An infinite IET $T$ is also measure preserving, being injective with
derivative~1 on~$I \setminus \Disc$, and $\Disc$ being of null measure.
Since $T$ is right-continuous on $I$, its values on $\Discaccr$ are determined
by its values on~$I\setminus \Discaccr$.
Like in the finite case, $\Disc$ can be replaced in the above definition by any
null-measure set $B$ of separation points such that $B \supset \Disc$, where $B
\cup \{1\}$ is compact.
Note that we can use a null-measure infinite set $B$ for finite
IETs too.

\subsection{Words} 

Let $A$ be a finite alphabet, let~$A^+$ be the set of finite words over~$A$,
let~$A^*$ be~$A^+$ with the empty word~$\epsilon$, and let~$A^\mathbb{N}$ be the
set of infinite words over~$A$.
The alphabet~$A$ can be endowed with the discrete topology, and the resulting
product topology on~$A^\mathbb{N}$ is metrizable and compact. The usual
associated metric is
the {\smbf Cantor metric}: if $w = a_0a_{1}a_{2}...$ and $w' = a'_0a'_1a'_2...$
in $A^\mathbb{N}$, with $a_i,a'_i \in A$, their distance is~0 if they are equal,
and~$2^{-k}$ if not, where $k$ is the smallest non-negative integer for which
$a_k \neq a'_k$.
Let~$v=v_0v_1...v_n \in A^*$, with $v_i \in A$, then its {\smbf cylinder set}
(for short, {\smbf cylinder}) is $Cyl(v)=\{w=a_0a_1a_2... \in A^\mathbb{N} \mid
a_0=v_0,a_1=v_1,...,a_n=v_n\}$.
A cylinder is a clopen set, and the collection of all the cylinders forms a
basis for the topology on~$A^\mathbb{N}$.  For~$\Lang \subseteq A^\mathbb{N}$
and $v \in A^*$, the {\smbf cylinder} in the subspace topology is $Cyl_\Lang(v)
= Cyl(v) \cap \Lang$.
Let~$\leq$ be an order over~$A$, propagated to all the words in $A^+$
and~$A^\mathbb{N}$ as the lexicographic order.  Then for every~$w,w' \in \Lang$
with $w \leq w'$, the {\smbf word interval} $[w, w']$ in $\Lang$ is $\{w'' \in
\Lang\:|\: w \leq w'' \leq w'\}$.
Note that the induced order topology is the same as the product topology.

A {\smbf factor} (or {\em subblock}) of a word~$w$ is a finite word $v$ such
that $w=w'vw''$, where $w', w''$ are possibly empty, and where $w'vw''$ denotes
the word concatenation of $w'$, $v$ and $w''$. The length of $v$ is denoted by
$|v|$.
The set of all distinct factors of a word~$w$ is denoted by~$Fact_w$, and for a
set $L \subseteq A^\mathbb{N}$, by $Fact_L=\bigcup_{w \in L} Fact_w$. For each
$n \in \mathbb{N}^*$, the set of all the distinct factors $v$ of~$w$ such that
$|v|=n$
is denoted by $Fact_w(n)$, and for $L \subseteq A^\mathbb{N}$, by
$Fact_L(n)=\bigcup_{w \in L} Fact_w(n)$.
The {\smbf (factor) complexity}~\cite{MH38,CN10} of a word~$w$ is the function
$p_w:\mathbb{N}^* \rightarrow \mathbb{N}^*$, with $p_w(n) = |Fact_w(n)|$, where
$|.|$ for a set denotes its cardinality, and for $L \subseteq A^\mathbb{N}$, it
is $p_L(n)=|Fact_L(n)|$. 
When $L$ is a shift, its {\smbf topological entropy}~\cite{Par66} is defined as
$\lim_{n\rightarrow \infty} log(p_L(n))/n$ (which exists, since $log(p_L(n))$
is subadditive).

An infinite word~$w$ is called {\smbf minimal}
when each factor in~$Fact_w$ occurs infinitely often in $w$ with bounded gaps,
i.e.~for each factor~$v$ of~$w=a_0a_1...$, the ordered sequence of distinct
indexes $\{n_j\}_{j \in \mathbb{N}}$ such that $a_{n_j}...a_{n_j+|v|-1} = v$ is
infinite, and $n_{j+1}-n_j$ is bounded independently of~$j$~\cite{MH38}.
Such a word $w$ is {\smbf minimal aperiodic} when it is not periodic, i.e.~when
there is no factor $v$ such that $w=v^\omega$.
A set of finite words is {\smbf prolongable} if all its words can always be
concatenated with letters to the right and to the left so that the resulting
words still belong to the set. When~$L \subseteq A^\mathbb{N}$ is made of
minimal words, $Fact_L$ is prolongable.
\subsection{Dynamical Systems and Shifts} 

Considering a self-map $f: X \rightarrow X$, where $X$ is a measured space and
where~$f$ is measurable, the pair $(X,f)$ is called a {\smbf measured dynamical
  system}. When $X$ is a topological space and $f$ is continuous, $(X, f)$ is
called a {\smbf topological dynamical system}.
The {\smbf shift map}~$\sigma$ over~$A^\mathbb{N}$ sends $a_0a_1...$ to
$a'_0a'_1...$, where $a'_i = a_{i+1}$ for every~$i \in \mathbb{N}$, and it is a
continuous map over~$A^\mathbb{N}$.  A {\smbf shift space} (or simply a {\smbf
  shift})~\cite{LM95} is a closed $\sigma$-invariant set of words $\Lang$
in~$A^\mathbb{N}$,
and accordingly, $(L,\sigma)$ is an instance of a topological dynamical system.
For a dynamical system $(X, f)$, the {\smbf (positive) orbit} of a point~$x \in
X$ is $\{f^n(x)\}_{n\in \mathbb{N}}$.
When $f$ is invertible, the {\smbf full~orbit} is $\{f^n(x)\}_{n\in
  \mathbb{Z}}$, and when $f$ is not invertible, it is the $\sim$-class of points
in~$X$ containing~$x$, where $x' \sim x$ if there are $n_1, n_2\geq 0$ such that
$f^{n_1}(x') = f^{n_2}(x)$~\cite{KST01}.
A continuous map $f: X \rightarrow X$ is called {\smbf minimal} if there is no
non-empty closed proper subset $X'\subsetneq X$ such that $f(X')=X'$, or
equivalently if the orbit of each $x\in X$ is dense in~$X$.
The map $f$ is {\smbf minimal aperiodic} if $X$ does not consist of a single
periodic orbit.  
Accordingly, a {\smbf minimal shift} is a shift containing no proper shift, and
a {\smbf minimal aperiodic shift} is a shift which is not made of the orbit of
one single periodic word.
If a word~$w$ belongs to a minimal shift~$L$, it is minimal as a word, and 
$Fact_w = Fact_L$.

For two measured dynamical systems $(X, f)$ and $(Y, g)$, if there exists a
measurable isomorphism $\phi: X \rightarrow Y$ 
such that $\phi \circ f = g \circ \phi$, then $(X, f)$ and $(Y, g)$ are said to
be {\smbf conjugate} by $\phi$.
When $(X, f)$ and $(Y, g)$ are topological dynamical systems, and $\phi$ is a
continuous onto map (resp. a homeomorphism), $(X, f)$ and $(Y, g)$ are said
to be {\smbf topologically semi-conjugate} (resp. {\smbf topologically
  conjugate}) by~$\phi$.
\section{From Shifts To Interval Exchange Transformations}

\subsection{The Conjugacies}
\label{conjugs_sec}

We start here by presenting basic results about how to embed
ordered measured spaces coming from shifts into $[0,1]$.

For any dynamical system $(X, f)$, where $X$ is a compact space and $f$ is a
continuous self-map on~$X$, there exist Borel probability measures on~$X$
which are $f$-invariant~\cite{HK02}.
If $f$ is minimal aperiodic, any such measure is nonatomic and takes positive
values on open sets.
Borel measures on compact sets behave well with respect to approximating
measurable sets by open and/or closed sets since they are {\smbf regular},
i.e.~for any measurable set~$E$,
$\mu(E)=\sup \{\mu(E') \;|\; E' \mbox{~\em compact}$, $E' \subseteq E\}$ and
$\mu(E)=\inf \{\mu(E') \;|\; E' \mbox{~\em open}$, $E \subseteq E'\}$.
For a regular measure, to be nonatomic is equivalent to every singleton having 
measure~0. 

From now on, $\Lang$ denotes a shift lexicographically ordered by~$\leq$,
and endowed with a $\sigma$-invariant Borel probability measure $\mu$ which is
nonatomic (so that the measure of any single word in $\Lang$ is zero), positive
on cylinders, and regular.  If $\Lang$ is a minimal aperiodic shift, such
measures exist.

\begin{Rem}\label{no_isolated_w_rem}
There is no isolated word in $\Lang$. 
\end{Rem}
\PR If there was an isolated word $w$ in $\Lang$, there would exist a prefix $v$
of $w$ such that $Cyl_\Lang(v)$ contains only $w$, but then $\{w\}$ would have
positive measure.\EPR 

Using $\mu$ and the order $\leq$ over $\Lang$, and denoting the smallest word of
$\Lang$ by~$w_{\Lang,min}$,
we define:

\[\begin{array}{lllllll}\label{valuation_phi0_def} 
 \phi_\mu: \;&\Lang \;\;\;\; &\rightarrow \;\;\;\; & \overline{I}=[0,1]\\
             &  w & \mapsto &\mu([w_{\Lang,min}, w]).
\end{array}\]

\begin{Lem}\label{non-decreas_lem} 
$\phi_\mu$ is a monotonic non-decreasing continuous map.
\end{Lem}
\PR Since $\mu$ is a measure which takes its values in $[0,1]$, $\phi_\mu$ is
monotonic non-decreasing.  For continuity, the involved spaces being compact,
we just check that the image by~$\phi_\mu$ of every sequence in $\Lang$
converging to some~$w \in \Lang$ is a sequence converging to $\phi_\mu(w)$.
Let $\cal S$ be such a sequence, from which we extract two
subsequences $\{w^+_i\}$ and~$\{w^-_i\}$, respectively formed by
decreasing words greater than~$w$ and by increasing words smaller
than~$w$.
At least one of them is infinite, say $\{w^+_i\}$, and 
$[w_{\Lang,min},w] = \bigcap_i [w_{\Lang,min},w^+_i)$.
Since $\mu$ is regular, $\mu([w_{\Lang,min},w])= \phi_\mu(w) =
\inf_i\mu([w_{\Lang,min},w^+_i))$, 
and since $\mu$ is zero on single words, for all $i>0$,
$\mu([w_{\Lang,min},w^+_i]) = \phi_\mu(w^+_i)$, thus
$\inf_i\phi_\mu(w^+_i)=\phi_\mu(w)$.
Hence since~$\phi_\mu$ is non-decreasing and since the infimum of
$\{\phi_\mu(w^+_{i})\}$ is $\phi_\mu(w)$, $\{\phi_\mu(w^+_{i})\}$
converges to~$\phi_\mu(w)$, whatever the choice of $\{w^+_i\}$.
If $\{w^-_i\}$ is also infinite, similar arguments as for $\{w^+_i\}$
apply using suprema instead of infima, whence the result.\EPR

\begin{Lem}\label{phioonto_lem} 
$\phi_\mu$ is a surjective map.
\end{Lem}
\PR Since $\Lang$ is compact and $\phi_\mu$ is continuous, 
$\overline{I} \setminus \phi_\mu(\Lang)$ is open, and if not empty it is a
disjoint union of open~intervals, since~0 and~1 belong to $\phi_\mu(\Lang)$.
Let $J=(x,x')$, with $x,x' \in \phi_\mu(\Lang)$, be one of these intervals, and
let $w_{x,sup} = \sup \{ w \in L \; | \; \phi_\mu(w)=x \}$ and $w_{x',inf} =
\inf \{ w \in L \; | \; \phi_\mu(w)=x' \}$.
Since $x<x'$ and since~$\phi_\mu$ is non-decreasing, $w_{x,sup} <
w_{x',inf}$. Also, there is no word $w \in \Lang$ with $w_{x,sup} < w <
w_{x',inf}$, otherwise $\phi_\mu(w_{x,sup}) \leq \phi_\mu(w) \leq
\phi_\mu(w_{x',inf})$, and by definition of $J$ we would have $\phi_\mu(w)= x$
or~$\phi_\mu(w)=x'$, contradicting the definitions of $w_{x,sup}$ or
$w_{x',inf}$.
Now, $\phi_\mu(w_{x',inf}) - \phi_\mu(w_{x,sup}) = \mu([w_{\Lang,min},
  w_{x',inf}]) - \mu([w_{\Lang,min}, w_{x,sup}]) = \mu((w_{x,sup},w_{x',inf}])$,
which is equal to the non-zero length of~$J$. But we just have checked
that~$(w_{x,sup},w_{x',inf}] = \{w_{x',inf}\}$ and~$\mu$ is zero on single
words. Hence there is no interval such as $J$, so $\phi_\mu$ is onto.  
\EPR

We now characterize the points where $\phi_\mu$ is non-injective. For
each~$n>0$, the length-$n$ factor set $Fact_\Lang(n)$ induces an ordered finite
partition of~$\Lang$ defined as

$$ CYL_{\Lang}(n): \;\;\;\;\Lang = \bigsqcup_{v \in Fact_{\Lang}(n)} Cyl_{\Lang}(v),
$$ 

\noindent
where the $Cyl_{\Lang}(v)$'s are ordered in~$\Lang$ according to the
lexicographic order of the~$v$'s in $Fact_{\Lang}(n)$.
By compactness, each~$Cyl_{\Lang}(v)$ has two {\smbf endpoints}: its smallest
and its greatest words.
We say that $w,w' \in \Lang$ (resp. $Cyl(u),Cyl(u') \in CYL_{\Lang}(n)$, $n>0$)
are {\smbf consecutive} if $w<w'$ (resp. $Cyl(u) < Cyl(u')$) and if there is no
word $w'' \in \Lang$ such that $w<w''<w'$ (resp. $Cyl(u) < w'' < Cyl(u')$).

\begin{Lem}\label{consec_words_lem}
Let $w,w' \in \Lang$. Then $w,w'$ are consecutive iff for some~$n > 0$ they are
endpoints, respectively the greatest and the smallest words, of two consecutive
cylinders of~$CYL_{\Lang}(n)$.
\end{Lem}
\PR 
~$(\Leftarrow)$: Trivial. 
~$(\Rightarrow)$: Assume $w < w'$, with $w,w'$ consecutive. Let $u \in A^*$ be
their longest common prefix, so that $w=ua_1a_2...$ and~$w'=ua'_1a'_2...$,
with~$a_i, a'_i \in A$.  Thus $w, w' \in Cyl_\Lang(u)$, and $a_1< a'_1$.
Since there is no word between~$w$ and~$w'$, the prefix~$ua_1...a_{j-1}a_j$ of
$w$ for each~$j > 1$, is also the greatest prefix for all the words in 
$Cyl_\Lang(ua_1)$.
Otherwise, there would be a factor~$ua_1v \in Fact_\Lang(|u|+j)$, $v\in A^+$,
such that $ua_1...a_{j-1}a_j < ua_1v$ as a prefix of a word~$w'' = ua_1v... \in
\Lang$. But then~$w < w'' < w'$, which is impossible.
Thus~$w$ is the greatest word in $Cyl_\Lang(ua_1)$. By similar arguments,
$w'$ is the smallest word in $Cyl_\Lang(ua'_1)$.
Finally, there is no other $a\in A$ with $a_1< a < a'_1$, such that there is
$w'''=ua... \in \Lang$, otherwise again $w < w''' < w'$.
Hence $Cyl_\Lang(ua_1)$ and $Cyl_\Lang(ua'_1)$ are consecutive cylinders in
$CYL_\Lang(|u|+1)$.
\EPR

\begin{Lem}\label{non-injectiv_lem} 
Let $w,w' \in \Lang$. Then $\phi_\mu$ is non-injective on $\{w,w'\}$ iff~$w,w'$
are consecutive.
\end{Lem}
\PR~$(\Leftarrow)$: Assume $w < w'$, with $w,w'$ consecutive. Then we have
$\phi_\mu(w') - \phi_\mu(w) = \mu([w_{\Lang,min}, w']) - \mu([w_{\Lang,min}, w])
= \mu((w, w']) = \mu(\{w'\}) = 0$,
since~$\mu$ is zero on single words. 
~$(\Rightarrow)$: 
Assume $w < w'$, with $w,w'$ non-consecutive. Thus there exists~$w'' \in \Lang$
with~$w < w'' < w'$, so that there is a length-$n$ prefix~$u$ of~$w''$, with
$n>0$, distinct from the length-$n$ prefixes of~$w$ and~$w'$.
But then, every~$w''' \in Cyl_\Lang(u)$ is such that~$w < w''' < w'$,   
and we have $\phi_\mu(w') - \phi_\mu(w) = \mu((w, w']) \geq \mu(Cyl_\Lang(u)) >
0$, since~$\mu$ is positive on cylinders.  
\EPR

Thus $\phi_\mu^{-1}(x)$, $x \in \overline{I}$, 
consists of either one or two words in~$L$, and a consequence of
Lemma~\ref{consec_words_lem} is that the set of points with two-words preimages
is countable.

Let us then transform~$\phi_\mu$ into an injective map by embedding
$\overline{I}$ as a subset of a larger compact space (using a similar
construction to the classic one for transforming piecewise continuous self-maps
into homeomorphisms, like e.g.~in~\cite{Kea75} for IETs). We first define the
following spaces and maps:

\begin{itemize}\itemsep=0cm 
\item $Z_0\subset \overline{I}$ denotes the image by~$\phi_\mu$ of the set of
  points where~$\phi_\mu$ is not injective.

\item $Z_0^-$ denotes a copy of~$Z_0$.

\item \label{Ipreiet_def} $\Ipreiet$ denotes $\overline{I} \sqcup Z_0^-$,
  ordered in such a way that each point in~$Z_0$ lies to the right of its copy
  in~$Z_0^-$, with no other point in between.  We endow~$\Ipreiet$ with the
  order topology for this order relation.
\item $\iota: \overline{I} \rightarrow \Ipreiet$ denotes the inclusion map. It is
  increasing and right-continuous on~$I$. 

\item $\kappa: \Ipreiet \rightarrow \overline{I}$ denotes the canonical map
  associated with the equivalence relation in~$\Ipreiet$ which identifies each
  point in~$Z_0$ with its copy in~$Z_0^-$. It is a non-decreasing continuous
  map, and it is onto. Accordingly, $\kappa \circ \iota$ is the identity map
  on~$\overline{I}$.
\end{itemize}

\noindent
Having in mind Lemma~\ref{non-injectiv_lem} we then define the following map
from $\phi_\mu$: {\small
\[\begin{array}{llll}
 \phi:\; & \Lang \;\;\;\; & \rightarrow \;\;\;\; & \Ipreiet \\ & \\
                    & w & \mapsto & \left\{
   \begin{array}{lll}
     \iota(\phi_\mu(w)) & \mbox{\em if $\phi_\mu(w) \notin Z_0$, or}\\ 
                      & \mbox{\em if $\phi_\mu(w) \in Z_0$, 
                   with $\phi_\mu^{-1}(\phi_\mu(w)) = \{w, w'\}$, $w>w'$.}\\
     \iota(\phi_\mu(w))^- & \mbox{\em if $\phi_\mu(w) \in Z_0$, 
                   with $\phi_\mu^{-1}(\phi_\mu(w)) = \{w, w'\}$, $w<w'$,} \\
                  & \mbox{\em where $\iota(\phi_\mu(w))^-$ is the copy
                              of $\iota(\phi_\mu(w))$ in $Z_0^-$.}
   \end{array} \right.
\end{array}\]
}

\noindent
Note that $\kappa \circ \phi= \phi_\mu$.

\vbox{
\begin{Lem}\label{phi_homeo_lem}
$\phi$ is an increasing homeomorphism. 
\end{Lem}
\PR By Lemma~\ref{non-decreas_lem}, $\phi$ is increasing where $\phi_\mu$ is
injective since $\phi = \phi_\mu$ on these points. 
Where $\phi_\mu$ is not injective, that is, for each $x \in Z_0$ where
$\phi_\mu^{-1}(x)=\{w,w'\}$ for some $w,w' \in \Lang$ with $w < w'$, then by
definition of~$\Ipreiet$ and $\phi$ we have $\phi(w) < \phi(w')$. Therefore
$\phi$ is everywhere increasing, thus into. Since $\Ipreiet =
\iota(\overline{I}) \sqcup Z_0^-$, $\phi$ is also onto. Hence $\phi$ is a
monotonic bijection between two totally ordered sets endowed with their
respective order topologies, so it is a homeomorphism.\EPR
}

\noindent
We define two more maps: 
{\small
\[\begin{array}{llllrllll}
 \Tpreiet:\; & \Ipreiet \;\;\;\; & \rightarrow \;\;\;\; & \Ipreiet &
 \;\;\;\;\;\;\;\;\;\;\;\; 
  \Tiet:\; & \overline{I} \;\;\;\; & \rightarrow \;\;\;\; & \overline{I}\\
    &x & \mapsto & \phi(\sigma({\phi}^{-1}(x))) &
    &x & \mapsto & \kappa(\Tpreiet(\iota(x))) 
\end{array}\]
}

\begin{Prop}\label{conjug_prop} 
Consider the following diagram:

\begin{center}
\begin{tikzpicture}
  \matrix (m) [matrix of math nodes,row sep=2.5em,column sep=3.5em,minimum width=2em]
              {
                L & \Ipreiet & \overline{I} \\
                L & \Ipreiet & \overline{I} \\
              };
      \draw[->] (m-1-1) edge node [left] {$\sigma$} (m-2-1);
      \draw[->] (m-1-2) edge node [left] {$\Tpreiet$} (m-2-2);
      \draw[->] (m-1-3) edge node [right] {$\Tiet$} (m-2-3);
      \draw[->] (m-1-1) edge node [above] {$\phi$} (m-1-2);
      \draw[->] (m-1-2) edge node [above] {$\kappa$} (m-1-3);
      \draw[->] (m-2-1) edge node [above] {$\phi$} (m-2-2);
      \draw[->] (m-2-2) edge node [above] {$\kappa$} (m-2-3);
      \draw[->] (m-1-1) edge [bend left=25] node [above] {$\phi_\mu$} (m-1-3);
      \draw[->] (m-2-1) edge [bend left=-20] node [below] {$\phi_\mu$} (m-2-3);
      \draw[->] (m-1-3) edge [bend left=20] node [below] {$\iota$} (m-1-2);
      \draw[->] (m-2-3) edge [bend left=-30] node [above] {$\iota$} (m-2-2);
\end{tikzpicture}
\end{center}

  \begin{enumerate}[a.]\itemsep=0cm
   \item For every $x \in I$, $\iota(\Tiet (x)) = \Tpreiet(\iota (x))$.
   \item $(L,\sigma)$ is topologically conjugate to $(\Ipreiet,\Tpreiet)$ by
     $\phi$.
   \item $((\iota(\overline{I}) \subset \Ipreiet),\Tpreiet)$ is topologically
     semi-conjugate to $(\overline{I},\Tiet)$ by $\kappa$.
   \item $((\phi^{-1}(\iota(\overline{I})) \subset L),\sigma)$ is topologically
     semi-conjugate to $(\overline{I},\Tiet)$ by $\phi_\mu$.
  \end{enumerate}
\end{Prop}
\PR (a): Note that $0, 1 \notin Z_0$.
Indeed, $\phi_\mu$ being non-decreasing, if $\phi_\mu^{-1}(0)$
(resp.~$\phi_\mu^{-1}(1)$) was made of two words, consecutive by
Lemma~\ref{non-injectiv_lem}, they would be the two smallest (resp.~greatest)
in~$\Lang$, and by Lemma~\ref{consec_words_lem} the smaller one (resp.~greater
one) would be isolated, contradicting Remark~\ref{no_isolated_w_rem}.
Now, we put $I'=\overline{I} \setminus (Z_0 \cup \{0, 1\})$. 
By definition, $\kappa$ is injective on $\iota(I)$, hence on $\iota(I')$,
thus since~$\kappa$ is a left inverse of $\iota$, it is also a right inverse
of~$\iota$ when restricted to this set.
We compose $\Tiet=\kappa \circ \Tpreiet \circ \iota$ on both sides to the left
by $\iota$, so that $\iota \circ \Tiet = \Tpreiet \circ \iota$ is valid on~$I'$,
provided that $\Tpreiet(\iota(I')) \subset \iota(I')$, 
which is proved as follows:
by definition of~$I'$, for every $x \in I'$,
$\phi_\mu^{-1}(x) = \phi^{-1}(\iota(x))$, that is, $\phi_\mu^{-1}(I') =
\phi^{-1}(\iota(I'))$.
Also, by Lemmas~\ref{consec_words_lem} and~\ref{non-injectiv_lem}, a word $w$
is in $\phi_\mu^{-1}(I')$ iff it is not a cylinder endpoint.
Thus, there are words in~$\Lang$ arbitrarily close to~$w$ from the right and the
left, that is, words smaller and greater than $w$ with arbitrarily long prefixes
common with the $w$'s ones.
By continuity of $\sigma$, the same holds for $\sigma(w)$, so 
$\sigma(\phi^{-1}(\iota(I'))) \subset \phi^{-1}(\iota(I'))$. 
Composing on both sides to the left by $\phi$ yields 
$\Tpreiet(\iota(I')) \subset \iota(I')$, as required.
Now, $Z_0$ being countable, $I'$ is dense in $\overline{I}$. Since both $\iota
\circ \Tiet$ and $\Tpreiet \circ \iota$ are compositions of right-continuous
maps on $I$, they are right-continuous too, so $\iota \circ \Tiet = \Tpreiet
\circ \iota$ extends to all of $I$.
\noindent
(b): 
$\phi$ is a homeomorphism, thus by definition of $\Tpreiet$, the result follows.
\noindent
(c): By definition of $\Tiet$, we have the conjugacy $\Tiet \circ \kappa =
\kappa \circ \Tpreiet$ on $\iota(I')$.
Both members of this identity are right-continuous maps, and by the same
arguments as in~(a), this conjugacy holds on all of~$\iota(I)$.  
Next, again by definition of $\Tiet$, we have $\Tiet(1) =
\kappa(\Tpreiet(\iota(1)))$, and $1$ has only one preimage by $\kappa$ which is
$\iota(1)$,
so the conjugacy holds on $\iota(\overline{I})$ too. Finally, $\kappa$ is
measure-preserving and onto.

\noindent
(d): Using the fact that $\kappa \circ \phi = \phi_\mu$, and the commutativity
properties given by~(b) and~(c) on $\kappa$ and $\phi$, the result follows.
\EPR

\subsection{Checking the Infinite Exchange Transformation Properties}

We now focus on when $\Tiet$ 
is an IET. 
Here are first three general lemmas: 

\begin{Lem}\label{T_preserves_Leb_lem}
The image of $\mu$ by $\phi_\mu$ is the Lebesgue measure on $\overline{I}$, which
is preserved by $\Tiet$.
\end{Lem}
\PR 
We denote by $f_* \mu$ the image of the measure $\mu$ by the map $f$.  By
Lemmas~\ref{non-decreas_lem} and~\ref{phioonto_lem}, $\phi_\mu$~is continuous and
surjective, thus the preimage of any interval $[0, x]$, with $x \in
\overline{I}$, is an interval $[w_{\Lang,min}, w]$ in $\Lang$, where $\phi_\mu(w)
= x = \mu([w_{\Lang,min}, w])$. Hence $({\phi_\mu})_* \mu$ is the Lebesgue measure
on $\overline{I}$.
Next, by Proposition~\ref{conjug_prop}(d), $\phi_\mu \circ \sigma = \Tiet \circ
\phi_\mu$ on $\phi^{-1}(\iota(\overline{I}))$, a set of full measure in
$\Lang$. Thus $({\phi_\mu})_* \sigma_* \mu = {\Tiet}_* ({\phi_\mu})_* \mu$,
and since~$\sigma$ preserves~$\mu$, $(\phi_\mu)_* \sigma_* \mu = ({\phi_\mu})_*
\mu$, that is,
$\Tiet$ preserves the Lebesgue measure on $\overline{I}$.  
\EPR

From $CYL_\Lang(1)$, that is, the partition of $\Lang$ defined as $\bigsqcup_{a
  \in A} Cyl_{\Lang}(a)$, we define the induced partition of $I$ as

\vspace*{-0.2cm}
\EQ P_A: \;\;\;\;I = \bigsqcup_{a \in A} \Iaiet, 
  \;\;\;\mbox{where } \Iaiet = \{x\in I\;|\: \phi^{-1}(\iota(x)) \in Cyl_\Lang(a)\}.
\label{part_PA_def}
\EEQ
\vspace*{-0.2cm}

\noindent
Note that each $\Iaiet$ is a right-open interval in $I$. Indeed, if $x < x'$,
with $x,x' \in \Iaiet$, then every $x''\in I$, with $x < x'' < x'$, is also
in~$\Iaiet$ since $\phi^{-1}$ and $\iota$ both preserve the order. Moreover, if
$a<a'$, where $a,a' \in A$ are consecutive, then if $w$ is the greatest word in
$Cyl_\Lang(a)$, and $w'$ the smallest in $Cyl_\Lang(a')$, we have
$\phi_\mu(w)=\phi_\mu(w')=y$, where $y\in Z_0$. By definition of $Z_0$ and
$\phi$ and since $w<w'$, $\phi^{-1}(\iota(y)) = w'$, thus $y \in I_{\Lang,a'}$,
which means that $\Iaiet$ is right-open.

\begin{Lem}\label{rightcontinuous_lem}
  $\Tiet$ is right-continuous on $I$, and piecewise increasing on $P_A$.
\end{Lem}
\PR $\Tiet$ is right-continuous on $I$, being a composition of the
right-continuous maps $\kappa, \Tpreiet$ and $\iota$.
Next, $\Tiet$ is non-decreasing on each $\Iaiet$ in $P_A$. Indeed, by
Proposition~\ref{conjug_prop}(d), $\Tiet\circ \phi_\mu=\phi_\mu\circ \sigma$ on
$\phi^{-1}(\iota(\Iaiet)) \subset Cyl_\Lang(a)$,
and $\phi_\mu$ is non-decreasing on~$I$, while $\sigma$ is increasing on
$Cyl_\Lang(a)$.
Now, assume there exist~$x,x'\in \Iaiet$, with $x<x'$, such that $\Tiet(x) =
\Tiet(x')$.  Since $\Tiet$ is non-decreasing it must be constant on $(x,x')$,
which is impossible since $\Tiet$ is measure-preserving
by~Lemma~\ref{T_preserves_Leb_lem}, hence $\Tiet$ is increasing on each $\Iaiet$.
\EPR

\begin{Lem}\label{selfmap_lem}
$\Tiet(I) \subset I$.
\end{Lem}
\PR By Lemma~\ref{rightcontinuous_lem}, $\Tiet$ is piecewise increasing on the
$\Iaiet$'s. Let~$J$ be one of these intervals, and assume $\Tiet(x)=1$ for some
$x \in J$. But $\Tiet$ being increasing on $J$, and $J$ being right-open,
$\Tiet(y)$ would be $>1$ for all $y\in J$ with $y>x$.~\EPR

As a next step, let us exhibit conditions so that there exists a partition
of~$I$ on which $\Tiet$ is a translation on its components ($P_A$ is not such a
partition in general).
For that purpose, we first study the effect of~$\sigma$ over the word intervals
in~$\Lang$.
A factor $u \in Fact_\Lang$ is called {\smbf left special} (resp.~{\smbf right
  special}) if $u$ has at least two distinct left (resp.~right) letter
prolongations in~$Fact_L$, i.e., if $u \in Fact_L(n)$, then $u$ is the suffix
(resp.~prefix) of at least two distinct factors in $Fact_L(n+1)$.

\begin{Lem}\label{shift_interval_lem} 
Let $\Lang$ be a shift such that $Fact_\Lang$ is prolongable.  Let $w,w' \in
Cyl_\Lang(au)$, with~$a \in A$,~$u \in Fact_L$ not left special, and $w <
w'$. Then

\vspace*{-0.3cm}
\begin{equation}
      \sigma ([w, w']) = [\sigma (w), \sigma(w')]. 
      \label{stab_sigma}
\end{equation}
\vspace*{-0.7cm}
\end{Lem}
\PR 
($\subseteq$): Let $w=auv$ and $w'=auv'$, with $v,v' \in A^\mathbb{N}$.  Let
$w''=auv'' \in [w, w']$, with $v'' \in A^\mathbb{N}$.
By lexicographic order, we have~$uv \leq uv''\leq uv'$.
Hence~$\sigma(w'')=uv''$ belongs to $[uv, uv']=[\sigma(w), \sigma(w')]$.

\noindent
($\supseteq$): Let again $w=auv$ and $w'=auv'$, with $v,v' \in
A^\mathbb{N}$. Let~$w''= uv'' \in [\sigma(w), \sigma(w')]$, with~$v'' \in
A^\mathbb{N}$.
Consider~$\{uu_j\}_{j\in\mathbb{N}}$ the sequence of all the prefixes of~$w''$
with~$u_j\in A^*$ starting with $u$.
By prolongability to the left, each~$uu_j$ has at least one left letter
prolongation, but it can only be~$a$ since~$u$ is not left special.
By prolongability to the right, each~$auu_j$ is the prefix of at least one word
in~$\Lang$, which makes a sequence converging to the word~$aw''$, belonging
to~$\Lang$ by compactness of~$\Lang$.  Thus there is one left letter
prolongation for~$w''$ in~$\Lang$,~unique since~$u$ is not left special, that
is, $aw''=\sigma^{-1}(w'')= \sigma^{-1}(uv'')$.
Since 
$uv \leq uv'' \leq uv'$, we have $auv \leq aw'' \leq auv'$, that is, $aw'' \in
[w, w']$. Hence $w'' \in \sigma ([w, w'])$.  \EPR

\noindent
Note that a specific case of the above is $\sigma(Cyl_\Lang(au)) =
Cyl_\Lang(u)$, where $a \in A$ and~$u$ is not left special.
Now, the idea is to build a partition of~$\Lang$ essentially made of cylinders
such that Equality~(\ref{stab_sigma}) holds in all of these cylinders.
An {\smbf infinite left special word} (or {\em infinite left special
  branch})
with respect to $Fact_L$, where $L\subseteq A^\mathbb{N}$, is a one-way infinite
word such that all its prefixes are left special factors in~$Fact_L$. We denote
by $\SP \subseteq A^\mathbb{N}$ the set of all the infinite left special words
with respect to~$Fact_L$.
When $\Lang$ is a shift, $\SP$ is included in $\Lang$.

\begin{Lem}\label{fundlemma}
Let~$\Lang$ be a shift having left special factors of arbitrary length,
measured by $\mu$, with $\mu(\SP)=0$. Then, there is an infinite
partition of $\Lang$ defined as
\vspace*{-0.1cm}
  \EQ
   PART_\Lang: \;\;\;\; \Lang = \bigsqcup_{k >0} Cyl_\Lang(v^{(k)}) \sqcup W_{\SP}, 
   \label{partL_def}
  \EEQ
\vspace*{-0.7cm}

\noindent
where:
\begin{itemize}\itemsep=0cm
\item For each $k>0$, $v^{(k)} \in Fact_\Lang$, such that for all $w,w' \in
  Cyl_\Lang(v^{(k)})$ with $w<w'$, $\sigma ([w, w']) = [\sigma(w), \sigma(w')]$
  (i.e.~Equality~(\ref{stab_sigma}) holds).
\item $W_{\SP}=\sigma^{-1}(\SP)$ is a null measure non-empty set in $\Lang$,
  which is formed by the accumulation of the endpoints of $Cyl_\Lang(v^{(k)})$
  in $\Lang$.
\end{itemize}
\end{Lem}
\PR 
Let us describe the cylinders of $PART_\Lang$ by an iterative process: as a
first step, consider the partition $CYL_{\Lang}(2)$ of~$\Lang$,
i.e.~$\bigsqcup_{a_ia_j \in Fact_\Lang(2)} Cyl_{\Lang}(a_ia_j)$, with $a_i,a_j
\in A$. By Lemma~\ref{shift_interval_lem}, for each non-left special~$a_j$, we
have $\sigma ([w, w']) = [\sigma (w), \sigma(w')]$ for all $w,w' \in
Cyl_\Lang(a_ia_j)$, so that $Cyl_\Lang(a_ia_j)$ is put in $PART_\Lang$, and
$a_ia_j$ is one of the~$v^{(k)}$'s.
As a second step, each of the remaining $Cyl_\Lang(a_ia_j)$'s not put in
$PART_\Lang$ during the first step, is partitioned with cylinders of the form
$Cyl_{\Lang}(a_ia_ja_k)$ in~$CYL_{\Lang}(3)$.
Again, for each non-left special suffix~$a_ja_k$, the word intervals in
$Cyl_\Lang(a_ia_ja_k)$ satisfy Equality~(\ref{stab_sigma}), and
$Cyl_\Lang(a_ia_ja_k)$ is put in $PART_\Lang$, while $a_ia_ja_k$ is also
one of the~$v^{(k)}$'s.
This refinement process is inductively applied as long as cylinders remain at
step $n$ by partitioning them with cylinders in $CYL_{\Lang}(n+1)$,
defining then all the cylinders of $PART_\Lang$.

Since the left special factors of $\Lang$ can have arbitrary length, and 
since every prefix of a left special word is left special, there is at least one
left special factor in~$\Lang$ for each length. Thus, given any $n>0$, and a
left special factor $u \in Fact_\Lang(n)$, the cylinder $Cyl_\Lang(au)$, $a
\in A$, must still be refined during the $n$th step of the above refinement
process. 
Therefore, this process is necessarily infinite. 
It determines infinite sequences of nested cylinders of the form
$\{Cyl_\Lang(au_j)\}_{j \in \mathbb{N}^*}$, where $u_j \in Fact_L(j)$ is left
special while being a prefix of~$u_{j+1}$, for all $j$.
The sequence $\{au_j\}_{j\in \mathbb{N}^*}$ gives a limit infinite word, which
belongs to $W_{\SP}$ since $\{u_j\}_{j\in \mathbb{N}^*}$ gives a limit infinite
word in $\SP$.
As a result, if a word $w \in \Lang$ does not belong to some
$Cyl_\Lang(v^{(k)})$, it belongs to $W_{\SP}$.
Conversely, if~$w \in \SP$, there exist more than one left letter prolongation
to each of its prefixes.  Consider any of these letters, say $a\in A$, so that
$aw \in \sigma^{-1}(\SP)=W_{\SP}$. According to the same above argument, each
prefix $u_j$ of~$w$ determines a cylinder $Cyl_\Lang(au_j)$ which has to be
refined further. We have $\bigcap_{\mbox{\footnotesize \em \;$u_j$ prefix of
    $w$}} Cyl_\Lang(au_j) = aw$, where $aw$ does not belong to any cylinder
$Cyl_\Lang(v^{(k)})$, and induces an infinite refinement.
Thus $PART_\Lang$ is a partition since a word $w$ in $\Lang$ belongs either to
one of the cylinders $Cyl_\Lang(v^{(k)})$ obtained through a finite number of
steps in the refinement process, or to~$W_{\SP}$ if $w$ induces an infinite
number of steps.

The set $W_{\SP}$ is not empty since $\SP$ is not empty when there are
infinitely many left special factors.
Also, $W_{\SP}$ is of null measure since
$\SP$ has been assumed of null measure and $\sigma$ is measure-preserving.
Finally, $W_{\SP}$ occurs as the set of the accumulation of the
$Cyl_\Lang(v^{(k)})$'s endpoints, since again~$\SP$ is of null measure, thus
there is no cylinder having all its words in $\SP$.
\EPR

\noindent
Note that if $\Lang$ is minimal aperiodic, it has left special factors of
arbitrary length.

In order to exhibit examples of aperiodic minimal shifts, we use the
following classic technique. 
A {\smbf substitution} is a map $\theta: A \rightarrow A'^*$, where $A$ and~$A'$
are alphabets, which is extended to words by sending
$w=... a_ia_{i+1}a_{i+2}...$ to
$\theta(w)=...\theta(a_i)\theta(a_{i+1})\theta(a_{i+2})...$, that is, on finite
words, $\theta$ is just a monoid morphism from~$A^*$ to~$A'^*$. For instance,
the {\smbf Thue-Morse substitution} $\theta_{tm}$ over $A=A'=\{a,b\}$ is defined
by $\theta_{tm}(a) = ab$ and~$\theta_{tm}(b) = ba$.
When $A=A'$, a substitution can be iterated.  The set of factors of such a
substitution is defined as $Fact_\theta =\{v \in A^*\;|\; v \in
Fact_{\theta^n(a)}, a \in A, n \in \mathbb{N}\}$, and
its {\smbf associated shift} as $\Lang_\theta=\{w \in A^\mathbb{N}\;|\; Fact_w
\subseteq Fact_\theta\}$.
When~$\theta$ is {\smbf non-erasing}, i.e.~there is no $a \in A$ whose image is
the empty word, and when $a$ is a strict prefix of $\theta(a)$, then for all
$n\geq 0$, $\theta^{n}(a)$ is a strict prefix of $\theta^{(n+1)}(a)$. Thus
$\{\theta^n(a)\}_{n \in \mathbb{N}}$ converges to a limit word in
$A^\mathbb{N}$, denoted by $\theta^\omega(a)$, which is a fixed point
of~$\theta$.
For instance, the Thue-Morse substitution has two fixed points:
$\theta_{tm}^\omega(a)=w_1= abbabaabbaababba...$, and
$\theta_{tm}^\omega(b)=w_2= baababbaabbabaab...$.
Considering a fixed point~$w$ of a substitution~$\theta$, its {\smbf associated
  shift} $L_{\theta,w}$ is then defined as the closure of $\{\sigma^n(w) \;|\; n
\in \mathbb{N}\}$.
A substitution~$\theta$ is said to be {\smbf primitive}\label{primitive_def} if
there~exists $n > 0$ such that for every $a,b \in A$, the word~$\theta^n(a)$
contains~$b$.
For a primitive substitution~$\theta$, it is known that~\cite{Que10}: (1)~if $w$
is any of its fixed points, $L_\theta= L_{\theta,w}$;
(2)~for all $n>0$, $L_{\theta^n}=L_\theta$, so that periodic points can also
be considered; (3)~$L_\theta$ is minimal.
Thus for instance the Thue-Morse substitution $\theta_{tm}$ generates a minimal
shift $L_{tm}$, equal to
$\Lang_{\theta_{tm}}=\Lang_{\theta_{tm},w_1}=\Lang_{\theta_{tm},w_2}$.  A
primitive substitution is said to be {\smbf aperiodic} if its associated shift
is minimal aperiodic, e.g.~$\theta_{tm}$ is known to be 
aperiodic~\cite{Thu12}\footnote{More generally, aperiodicity for a substitution
  is decidable using~\cite{ER83}.}.

\begin{Example} (Step-by-step construction of a partition $PART_\Lang$).\em
\label{thue-morse-ex}
Considering the above Thue-Morse shift~$L_{tm}$, and its fixed points $w_1,
w_2$, we build here the first components of $PART_{\Lang_{tm}}$ as in the proof
of Lemma~\ref{fundlemma}. Its first factors are:

\vspace*{-0.6cm}
\[\begin{array}{lllllll} 
  Fact_{\Lang_{tm}}(1) &=&\{a,b\}\\
  Fact_{\Lang_{tm}}(2) &=&\{aa,ab,ba,bb\}\\
  Fact_{\Lang_{tm}}(3) &=&\{aab, aba, abb, baa, bab, bba\}\\
  Fact_{\Lang_{tm}}(4) &=&\{aaba, aabb, abaa, abab, abba, baab, baba, babb, bbaa\}\\
  Fact_{\Lang_{tm}}(5) &=&\{aabab, aabba, abaab, ababb, abbaa, abbab, baaba,
                            baabb, babaa,\\
                       & &\;\;babba, bbaab, bbaba\}
\end{array} \]
\vspace*{-0.4cm}

\noindent
Thus its first left special factors are:
$a, b, ab, ba, aba, abb, baa, bab, abba, baab.$
Applying the refinement process, the first step, based on $CYL_\Lang(2)$, gives
no cylinders of $PART_{\Lang_{tm}}$ since $a$ and $b$ are both left special
factors. The second step, based on $CYL_\Lang(3)$, gives $Cyl_{\Lang_{tm}}(abb)$
and $Cyl_{\Lang_{tm}}(baa)$, since $aa$ and $bb$ are not left special factors,
thus $v^{(1)}=abb$ and $v^{(2)}=baa$. The third step gives no new cylinders
since $abba$ and $baab$ are the only factors whose right suffix is not left
special, already contained in the cylinders of the preceding step.  The fourth
step gives $Cyl_{\Lang_{tm}}(aabab)$, $Cyl_{\Lang_{tm}}(ababb)$,
$Cyl_{\Lang_{tm}}(babaa)$, and $Cyl_{\Lang_{tm}}(bbaba)$, thus $v^{(3)}=aabab$,
$v^{(4)}=ababb$, $v^{(5)}=babaa$, and $v^{(6)}=bbaba$. And so on. \EEX

Now, given a shift $\Lang$, and $PART_\Lang$ as obtained from
Lemma~\ref{fundlemma}, we define the induced partition of~$I$ as follows:

\vspace*{-0.3cm}
\EQ
 PART_{\Lang,I}: \;\;\;\; I = \bigsqcup_{k \in \mathbb{N}^*} \Ikiet \sqcup B_{\SP},
 \label{partI_def}
\EEQ 
\vspace*{-0.2cm}

\noindent
where $\Ikiet = \{x \in I\;|\: \phi^{-1}(\iota(x)) \in Cyl_\Lang(v^{(k)}) \}$,
and $B_{\SP}=\{x \in I\;|\: \phi^{-1}(\iota(x)) \in W_{\SP}\}$.
Each $\Ikiet$ is a right-open interval by the same arguments as for the intervals
of the partition $P_A$ (see p.~\pageref{part_PA_def}).

Let~$\Biet$ be 
the closure in $\overline{I}$ of the left ends of the $\Ikiet$'s, intersected
with $I$ (so that $1 \notin \Biet$).
Let $\Bietacc$ (resp. $\Bietaccr$) be the set of accumulation points (resp.
accumulation points from the right) of $\Biet$ in~$\overline{I}$. 
Note that $B_{\SP} \subset \Bietacc$, and even that $B_{\SP}=\Bietaccr$ since an
accumulation point $x$ only to the left is the left endpoint of some $\Ikiet$
and does not belong to $B_{\SP}$ (in this case, $\phi^{-1}(\iota(x)) \notin
W_{\SP}$ because of the differences between~$\iota(I)$ and~$\Ipreiet$ --~see
also Definition~\ref{iiet_def}).

\begin{Lem}\label{separ_nullmeas_lem}
$\Biet \sqcup \{1\}$ is a null-measure infinite compact set of points in
  $\overline{I}$.
\end{Lem}
\PR The set $\Biet \sqcup \{1\}$ is compact by definition.  Since the set of the
left endpoints of the~$\Ikiet$'s is countable, the measure of $\Biet$ is the
same as that of $\Bietacc$, and also as that of $\Bietaccr$ which is equal to
$B_{\SP}$.
Now, recall from~Lemma~\ref{T_preserves_Leb_lem} that the Lebesgue measure is
the image of $\mu$ by $\phi_\mu$. Since $\phi_\mu^{-1}(B_{\SP}) = W_{\SP}$ and
$\mu(W_{\SP}) = 0$, the measure of $\Biet$ is also zero in $\overline{I}$.
\EPR

\begin{Lem}\label{veechlemma}
$\Tiet$ is a translation
on each interval component $\Ikiet$ of $PART_{\Lang,I}$.
\end{Lem}
\PR Let the interior of $\Ikiet$ be denoted by $I^{(k)^o}_\Lang$. By
definition of~$\Tiet$, for all $x \in I^{(k)^o}_\Lang$,

\vspace*{-0.7cm}
\EQ
\Tiet(x)-x =\kappa(\phi(\sigma({\phi}^{-1}(\iota(x))))) - x.
\label{diff_basic_eq}
\EEQ

\noindent
We first put ${\phi}^{-1}(\iota(x)) = w$.  Next, recalling that $\phi_\mu =
\kappa \circ \phi$, and that $\kappa \circ \iota$ is the identity, we compose
${\phi}^{-1}(\iota(x))=w$ on both sides to the left by $\kappa\circ\phi$, and
get $x=\kappa(\phi (w)) = \phi_\mu(w)$, so that~(\ref{diff_basic_eq}) becomes
for all $x \in I^{(k)^o}_\Lang$, 
\EQ \Tiet(x)-x = \kappa({\phi}(\sigma(w))) - \kappa({\phi}(w))  
       = {\phi_\mu}(\sigma(w)) - {\phi_\mu}(w).
\label{diff_basic2_eq}
\EEQ
\vspace*{-0.5cm}

\noindent
By definition of $PART_{\Lang,I}$, the inverse image of $I^{(k)^o}_\Lang$ by
$\phi_\mu$ is the interior $C^{(k)^o}$, where $C^{(k)}$ denotes
$Cyl_\Lang(v^{(k)})$ of $PART_\Lang$.
As a closed ordered set, $C^{(k)}$ has a smallest word $w_{C^{(k)},min}$,
and recall that $w_{\Lang,min}$ denotes the smallest word of~$\Lang$, so that
for all $w \in C^{(k)^o}$, $[w_{\Lang,min}, w] = [w_{\Lang,min},
  w_{C^{(k)},min}] \sqcup (w_{C^{(k)},min}, w]$, then 
\EQ
{\phi_\mu}(w) = \mu([w_{\Lang,min}, w_{C^{(k)},min}]) + 
                  \mu([w_{C^{(k)},min}, w]). 
\label{phi_o_eq} 
\EEQ

\noindent
Since $w_{C^{(k)},min}$ and $w$ belong to the same cylinder, and
$w_{C^{(k)},min}< w$,
we have $\sigma(w_{C^{(k)},min}) < \sigma(w)$, hence
$$ {\phi_\mu}(\sigma(w)) = \mu([w_{\Lang,min}, \sigma(w_{C^{(k)},min})]) + 
          \mu([\sigma(w_{C^{(k)},min}), \sigma(w)]).$$

\noindent
Next, using the properties over the components of $PART_\Lang$, 
\[\begin{array}{lll}\label{phi_o_sigma_eq} 
 {\phi_\mu}(\sigma(w)) &= 
           \mu([w_{\Lang,min}, \sigma(w_{C^{(k)},min})]) + 
           \mu(\sigma([w_{C^{(k)},min},w]))  
              \;\;\;\;\mbox{\small (by Equality~(\ref{stab_sigma}))}\\
  &=  \mu([w_{\Lang,min}, \sigma(w_{C^{(k)},min})]) + 
                          \mu(\sigma^{-1}\sigma([w_{C^{(k)},min},w]))
              \;\;\mbox{\small ($\mu$-preservation)}\\
  &= \mu([w_{\Lang,min}, \sigma(w_{C^{(k)},min})]) + 
                         \mu([w_{C^{(k)},min},w]).
              \;\;\;\;\mbox{\small ($\sigma$ injectivity on $C^{(k)^o}$)} \\
\end{array}\]

\vspace*{-0.1cm}
\noindent
Therefore, from the above last equation and~(\ref{phi_o_eq}), 
we get for all $w \in C^{(k)^o}$, 
\[\begin{array}{llll}
    {\phi_\mu}(\sigma(w)) - {\phi_\mu}(w) 
  & = &\mu([w_{\Lang,min}, \sigma(w_{C^{(k)},min})]) - 
                                  \mu([w_{\Lang,min}, w_{C^{(k)},min}]) =
     K_{C^{(k)}.}
\end{array}\]

\vbox{
\noindent
Thus $K_{C^{(k)}}$ depends only on $C^{(k)^o}$ to which $w$ belongs,
that~is, it depends only on the open interval $I^{(k)^o}_\Lang$ to which $x$
belongs. Hence, for every $\Ikiet$ of $PART_{\Lang,I}$, and its corresponding
cylinder $C^{(k)}=Cyl_\Lang(v^{(k)})$ of $PART_\Lang$, there is a
constant~$K_{C^{(k)}}$ such that for all $x \in I^{(k)^o}_\Lang$, $\Tiet(x)-x =
K_{C^{(k)}}$. By right-continuity of~$\Tiet$ (see
Lemma~\ref{rightcontinuous_lem}), this equality extends to all of $\Ikiet$.
\EPR
}

Following Definition~\ref{iiet_def}, let $\Disciet$ be the closure of
$Y_\Lang$ in $I$, where $Y_\Lang$ is the set of
discontinuities of~$\Tiet$, which by Lemma~\ref{veechlemma}, is necessarily
included into~$\Biet$.  Let $\Discietaccr$ be the set of accumulation points of
$\Disciet$ from the right.

\begin{Lem}\label{injectiv_disc_acc_lem}
$\Tiet$ is injective on $I \setminus \Discietaccr$.
\end{Lem}
\PR 
By Lemma~\ref{veechlemma}, $\Tiet$ is a translation on each of the right-open
interval components $\Ikiet$ of $PART_{\Lang,I}$.
The set $\Disciet$ can be obtained by dropping all the points from $\Biet$
having a neighborhood where $\Tiet$ extends as a translation map, so that
$\Tiet$ is still a translation map on each right-open interval component of~$I
\setminus \Disciet$.
By Lemma~\ref{T_preserves_Leb_lem}, the Lebesgue measure is preserved by
$\Tiet$, and by Lemma~\ref{separ_nullmeas_lem}, $\Biet$ has null measure, 
thus $\Disciet \subseteq \Biet$ has too.
Now, if the open interval components of $I \setminus \Disciet$ had overlaps
in their translated images by~$\Tiet$, 
the global image would be a set $I'$ with $\mu(I')<1$, such that
$\mu(T^{-1}(I')) = \mu(I \setminus \Disciet) = 1$, contradicting measure
preservation.
Moreover, no image by $\Tiet$ of a left endpoint of these open intervals can lie
in the image of another interval, thus $T$ is injective on all the corresponding
right-open intervals. Since $\Disciet \setminus \Discietaccr$ consists of all
these left endpoints, the result follows.
\EPR

\begin{Lem}\label{nullmeasure_infinspec_lem}
Let $L$ be a shift measured by $\mu$.  If $\Lang$ has zero topological entropy,
then $\mu(\SP)=0$.
\end{Lem}
\PR The topological entropy dominates the (measure-theoretic) entropies
on~$\Lang$ with respect to all the invariant probability
measures~\cite[Proposition 4.4.1]{HK02}. Thus with respect to~$\mu$, the shift
$\Lang$ has necessarily zero entropy.
But then, a dynamical system with zero entropy with respect to an invariant
measure is known to be invertible~\cite[Section 3.7k]{HK02}), that is, it has an
inverse on a subset of full measure. In the case of a shift $\Lang$, the set of
words on which the shift map is not invertible is exactly $\SP$, so that
$\mu(\SP)=0$.  \EPR

We are now ready to prove Theorem~\ref{main1} as stated in the introduction: 

\PRnotheo Since~$\Lang$ is assumed aperiodic minimal, the results 
of Section~\ref{conjugs_sec} apply,
so that $(\Lang, \sigma)$ is a topological conjugate by~$\phi$ of
$(\Ipreiet,\Tpreiet)$, where~$I$ is embedded as a subset of~$\Ipreiet$ in such a
way that~$\Tpreiet$ is an extension of~$\Tiet$. By
Lemma~\ref{rightcontinuous_lem}, $\Tiet$ is right-continuous, and by
Lemma~\ref{nullmeasure_infinspec_lem}, $\mu(\SP)=0$ for any $\sigma$-invariant
Borel probability measure $\mu$ chosen for $\Lang$. Minimal aperiodicity
means $\Lang$ has left special factors of arbitrary length, thus
Lemma~\ref{fundlemma} applies. Hence, by Lemmas~\ref{separ_nullmeas_lem},
\ref{veechlemma} and~\ref{injectiv_disc_acc_lem}, and since the set $\Disciet$
of $\Tiet$ is included into $\Biet$, $\Tiet$ is an IET.
\EPR 

Note that Theorem~\ref{main1} remains valid under the assumptions on~$\Lang$ in
Lemma~\ref{fundlemma}, in particular $\mu(\SP)=0$ can be used instead of zero
topological entropy.
Note also that the above construction relies on $\Biet$, a set determined from
$PART_{\Lang,I}$,
determined in turn from $PART_\Lang$.
This set of separation points for~$\Tiet$ could have been different, as it is
just a set with the property of including~$\Disciet$ (see
Section~\ref{iet_subsec}). However, following the construction of $PART_\Lang$
in the proof of Lemma~\ref{fundlemma}, this set has some word-combinatorics
properties, and we shall use them henceforth.

As a complement, in view of Definition~\ref{iiet_def} of an infinite IET and of
Lemma~\ref{injectiv_disc_acc_lem}, we give the following example:

\begin{Example} ($\Tiet$ can be non-injective on~$\Discietaccr$). \em
\label{non-injectivity-ex}
Let us reconsider the Thue-Morse substitution $\theta_{tm}$ over $A=\{a,b\}$
(see Example~\ref{thue-morse-ex}), its associated shift~$\Lang_{tm}$, and its
fixed points $w_1=\theta_{tm}^\omega(a)$ and $w_2=\theta_{tm}^\omega(b)$.
A {\smbf bispecial factor} is a factor which is both right and left special.
According to~\cite{Cas95}, the bispecial factors in $\Lang_{tm}$ are
$\theta_{tm}^n(a)$, $\theta_{tm}^n(b)$, $\theta_{tm}^n(aba)$,
$\theta_{tm}^n(bab)$, with $n \geq 0$.
Thus we have $\{w_1,w_2\} = SP_{L_{tm}}$. In fact, in an aperiodic minimal
shift, every infinite special word has always infinitely many prefixes which are
bispecial factors.

Let the lexicographic order over $\Lang_{tm}$ be determined by $a<b$.  We
claim that $w_1,w_2$ are consecutive in $\Lang_{tm}$, with $w_1$ being the
greatest word in $Cyl_{\Lang_{tm}}(a)$, and $w_2$ the smallest word in
$Cyl_{\Lang_{tm}}(b)$.
Assume on the contrary there is $w'_1 \in Cyl_{\Lang_{tm}}(a)$ with $w_1<w'_1$.
Let $u$ be the maximal common prefix so that $w_1=uv_1$ and $w'_1=uv'_1$. Then
$u$ is bispecial, since at least three out of the factors $aua$, $aub$, $bua$,
$bub$ are in $Fact_{\Lang_{tm}}$ (both $ua$ and $ub$ are prolongable to the
left, and at least one of them is left special, being a prefix of $w_1$). 
By~\cite{Cas95}, such a bispecial factor in $\Lang_{tm}$ with at least three
continuation factors, and starting by $a$ must be of the form
$\theta_{tm}^n(a)$, for some $n \geq 0$.
Now, since $w_1<w'_1$, the word~$v_1$ must also start with~$a$, and $v'_1$
with~$b$.
However since there is $n>0$ such that $u=\theta_{tm}^{n}(a)$,
and since $\theta_{tm}^{n+1}(a) = \theta_{tm}^{n}(a)\theta_{tm}^{n}(b)$ is
also a prefix of~$w_1$, then $u$ must be followed by~$b$ in $w_1$.
Hence there is no word as~$w'_1$, and~$w_1$ is the greatest word in
$Cyl_{\Lang_{tm}}(a)$.
Similar arguments apply to~$w_2$ showing it is the smallest word in
$Cyl_{\Lang_{tm}}(b)$.
Thus by Lemma~\ref{consec_words_lem}, $w_1$, $w_2$ are consecutive.

As a result, $\{aw_1,aw_2\}$ and $\{bw_1,bw_2\}$ are two pairs of consecutive
words in $\Lang_{tm}$.
Moreover, $Cyl_{\Lang_{tm}}(ab)$ 
is a non-empty cylinder between $aw_2$ and $bw_1$, then by
Lemmas~\ref{non-decreas_lem} and~\ref{non-injectiv_lem},
$\phi_\mu(aw_1)=\phi_\mu(aw_2)<\phi_\mu(bw_1)=\phi_\mu(bw_2)$.
Now, since $\Lang_{tm}$ has linear complexity~\cite{Cas95}, it has zero
topological entropy, and Theorem~\ref{main1} applies. We~then obtain a conjugate
IET $\Tiettm$ with $B_{L_{tm}}$ as set of separation points. By
Proposition~\ref{conjug_prop}(d), and since $\{aw_2, bw_2\} \subset
\phi^{-1}(\iota(\overline{I}))$ by definition of~$Z_0$,
we have $\Tiettm(\phi_\mu(aw_2)) = \phi_\mu(\sigma(aw_2)) =\phi_\mu(w_2)=
\phi_\mu(\sigma(bw_2))= \Tiettm(\phi_\mu(bw_2))$.
Hence $\Tiettm$ is not injective.

Now, by Lemma~\ref{injectiv_disc_acc_lem}, $\phi_\mu(aw_2)$ and $\phi_\mu(bw_2)$
must belong to the set of discontinuities ${\cal D}_{\Lang_{tm},acc,r}$.
Note that this is coherent with the fact that by Lemma~\ref{fundlemma},
$\{aw_1,bw_1,aw_2,bw_2\} = W_{SP_{\Lang_{tm}}}$,
so that $\{\phi_\mu(aw_2), \phi_\mu(bw_2)\} = B_{SP_{L_{tm}}}$, and that
by definition, $B_{SP_{L_{tm}}} = B_{\Lang_{tm},acc,r} \supseteq {\cal
  D}_{\Lang_{tm},acc,r}$.
\EEX

Also, according to the construction behind Lemma~\ref{fundlemma}, the set of
separation points~$\Biet$ is always infinite, while including the set $\Disciet$
which can be in fact finite (see e.g.
Corollary~\ref{finite_iet_gives_the_same_finite_iet_cor} further), but:

\begin{Example} (The set $\Disciet$ of $\Tiet$ can be infinite). \em
\label{infinite-discontinuities}
Let us consider again the Thue-Morse associated shift $L_{tm}$ and its
associated definitions.
If there were only finitely many discontinuities for $\Tiettm$, we could assume
there is an open interval $J \subset I$,
having $\phi_\mu(a w_2)$ as left endpoint, on which $\Tiettm$ is continuous.  We
could also assume~$J$ avoids the finite set $\phi_\mu(W_{SP_{L_{tm}}})$, hence~$J$
would be contained in a subset of~$I$ where~$\Tiettm$ is injective.
But $\Tiettm(\phi_\mu (a w_2)) = \Tiettm(\phi_\mu (b w_2))$, so by
right-continuity, $\Tiettm$ cannot be injective on~$J$, and there is no such
interval.
\EEX

\subsection{Codings and Reconstructions}
\label{coding_reconstruct_sec}

In this section, we describe some more relationships between $\Lang$ and $\Tiet$. 

First, an IET $T$ is said to be {\smbf minimal} if all its orbits are dense
in~$I$.
Next, a {\smbf monotonicity partition}~$P$ of~$I$ with respect to~$T$ is a
finite partition, made of right-open intervals on each of which $T$ is
increasing, but not necessarily continuous.
Recall then that by Lemma~\ref{rightcontinuous_lem}, the conjugate IETs given by
Theorem~\ref{main1} are piecewise increasing, thus admitting monotonicity
partitions.
The associated {\smbf coding} of the orbits of $T$ with a partition $P$ is based
on a map~$\alpha$ assigning a distinct letter of an alphabet~$A$ to each
interval of $P$, so that 

\vspace*{-0.2cm}
\[\begin{array}{llll} 
 cod_{0,P}:\; & I \;\;\;\; & \rightarrow \;\;\;\; & A \\
          & x          & \mapsto &\mbox{\em $\alpha(J)$ \;\;\;if $x \in J$, with
                                         $J \in P$,} 
\label{cod_0_P_def}
\end{array}\]

\noindent
which is extended as 

\vspace*{-0.2cm}
\[\begin{array}{lllllll} 
 cod_{P}:\; & I \;\;\;\; & \rightarrow \;\;\;\; & A^\mathbb{N} \\ & x & \mapsto &
 cod_{0,P}(x)\; cod_{0,P}(T(x))\; cod_{0,P}(T^2(x))...
\end{array}\]

\noindent
The word $cod_P(x)$ is called the {\smbf symbolic orbit} of $x$ by $cod_P$, and
the closure $\overline{cod_P(I)}$, denoted by $\Lang_{P}$, is called the {\smbf
  associated shift} of $cod_P$.
A~lexicographic order on $\Lang_P$ is induced by the order over $A$ determined
by the order of the component intervals of $P$ on~$I$, that is, $cod_{0,P}(x) <
cod_{0,P}(x')$, if $x\in J$, $x'\in J'$ with $J < J'$, i.e.~$J$ occurs before
$J'$ on $I$.

The next two technical lemmas make it possible to apply Theorem~\ref{main1}
to~$\Lang_P$: 

\begin{Lem}\label{cod_increasing_lem}
  Let $T$ be a piecewise increasing IET, and let $P$ be a monotonicity
  partition of~$I$. Then $cod_P$ is non-decreasing, Borel measurable and
  right-continuous.  Moreover, if $T$ is minimal, $cod_P$ is increasing.
\end{Lem}
\PR Let $x,x' \in I$, with $x<x'$. If $x \in J$ and $x'\in J$, where $J,J'$ are
intervals of $P$ with $J<J'$, we have $cod_{0,P}(x) < cod_{0,P}(x')$, hence
$cod_P(x) < cod_P(x')$. Otherwise, $J=J'$, so that $cod_{0,P}(x)=cod_{0,P}(x')$.
Since on the intervals of $P$, we have $T(x) < T(x')$, we can inductively apply
the same argument until for some~$n > 0$, $cod_{0,P}(T^n(x)) <
cod_{0,P}(T^n(x'))$, in which case $cod_P(x) < cod_P(x')$,
and if there is no such~$n$, then $cod_P(x) = cod_P(x')$. Hence, $cod_P$ is
non-decreasing.
If $T$ is minimal, the above~$n$ always exists since $T^n(x)$ can be
arbitrarily close to the left end of an interval in $P$, and as long as both
$T^n(x)$ and $T^n(x')$ are in the same interval of $P$, we have $T^n(x') -
T^n(x) \geq x' - x$. Hence, $cod_P(x) < cod_P(x')$.

The map $cod_P$ is the product of the maps $cod_{0,P} \circ T^i$, $i \ge 0$,
each being~continuous on the interiors of the intervals of a finite partition
of~$I$, since $T^i$ is piecewise increasing.
Hence $cod_P$ is Borel measurable.
Likewise, each $cod_{0,P} \circ T^i$ is right-continuous since $cod_P$ and $T^i$
are right-continuous, hence $cod_P$ is too.
\EPR

\noindent
From the above lemma, the image by $cod_P$ of the Lebesgue measure on 
$I$ is a Borel measure on $\Lang_P$, that we denote by $\mu_P$.
Since $cod_P$ conjugates $\sigma$ and $T$, $\mu_P$~is invariant by~$\sigma$.
Also, $cod_P$ being right-continuous, when~$T$ is minimal, $\Lang_P$ is
aperiodic minimal, and $\mu_P$ is nonatomic
and positive on the cylinder~sets of~$\Lang_P$.

\begin{Lem}\label{SP_measure_zero_lem}
  Let $T$ be a piecewise increasing and minimal IET.  Let $P$ be a monotonicity
  partition of $I$, and let $\Lang_{P}$ be the associated shift. Then
  $\mu_P(SP_{\Lang_{P}}) = 0$.
\end{Lem}
\PR By definition of $\mu_P$, for all $x \in I$,
$\mu_P([w_{\Lang_{P},min}, cod_P(x)])$ is equal to the Lebesgue measure of
$cod_P^{-1}([w_{\Lang_{P},min}, cod_P(x)))$,
that is, of $[0, x)$, since by Lemma~\ref{cod_increasing_lem} $cod_P$ is
  increasing.
Thus, using the map $\phi_{\mu_P}$ on $\Lang_{P}$, for all $i\in I$, 
\vspace*{-0.2cm}
\EQ
\phi_{\mu_P}(cod_P(x)) = x.
\label{pushf_basic_eq1}
\EEQ
\vspace*{-0.7cm}

\noindent
By Lemma~\ref{non-injectiv_lem} and by definition of $Z_0$, recall that if $x
\notin Z_0$, $\phi^{-1}_{\mu_P}(x)$ is well-defined, and by
Equality~(\ref{pushf_basic_eq1}) is equal to $cod_P(x)$. If $x\in Z_0$,
$\phi_{\mu_P}^{-1}(x)$ consists of two words, one of them being $cod_P(x)$.
Accordingly, $cod_P(I)$ misses some words of $\Lang_P$, but $\Lang_P \setminus
cod_P(I)$ is countable, since $Z_0$ is countable.

Let us now estimate the measure of $\sigma^{-1}(SP_{\Lang_P})$.
We define the subset $W = \sigma^{-1}(SP_{\Lang_P}) \cap (\Lang_P \setminus
cod_P(I))$,
which is countable too.
Let $w_0 \in \sigma^{-1}(SP_{\Lang_P})$. By definition of $SP_{\Lang_P}$,
$\{w_0\} \subsetneq \sigma^{-1}(\sigma(w_0))$,
and there are three cases:

\vspace*{-0.1cm}
\begin{itemize}\itemsep=-0.1cm
\item $w_0 \in W$: thus $w_0 \in \sigma^{-1}(\sigma(W))$.

\item $w_0 \notin W$ and there is $w_1 \in \sigma^{-1}(\sigma(w_0))$ such that
  $w_1 \in W$: thus $w_0 \in \sigma^{-1}(\sigma(W))$ too, since $\sigma(w_0) =
  \sigma(w_1)$.

\item $w_0 \notin W$ and there is no $w_1 \in \sigma^{-1}(\sigma(w_0))$ such
  that $w_1 \in W$: thus there are $x_0,x_1 \in I$
   such that $\{w_0, w_1\} \subset \sigma^{-1}(\sigma(w_0))$ with $w_0 =
   cod_P(x_0)$, $w_1 = cod_P(x_1)$, and $x_0 \neq x_1$ since $w_0 \neq w_1$.
  By definition of $cod_P$, for all $x \in I$, 

\vspace*{-0.4cm}
\EQ
  cod_P(T(x)) = \sigma(cod_P(x)). 
\label{conjug_basic_eq1}
\EEQ
\vspace*{-0.7cm}

\noindent
Since $\sigma(cod_P(x_0))=\sigma(cod_P(x_1))$, then $cod_P(T(x_0)) =
cod_P(T(x_1))$, that is, since $cod_P$ is injective, $T(x_0)=T(x_1)$.  Hence
according to Definition~\ref{iiet_def}, $x_0, x_1 \in \Discaccr$
 and $w_0 \in \sigma^{-1}(cod_P(\Discaccr))$. 
\end{itemize}

\noindent
Thus $\sigma^{-1}(SP_{\Lang_P}) \subseteq \sigma^{-1}(\sigma(W)) \cup
\sigma^{-1}(cod_P(\Discaccr))$.
Then, $\mu_P(\sigma^{-1}(\sigma(W))=0$, since $\sigma^{-1}(\sigma(W))$ is
countable
and $\mu_P$ is nonatomic.
Next, $\mu_P(\sigma^{-1}(cod_P(\Discaccr)))=\mu_P(cod_P(\Discaccr))$ since
$\mu_P$ is preserved by $\sigma$. Now, since $\mu_P$ is the image of the
Lebesgue measure on~$I$ by $cod_P$ and $\Discaccr$ has measure~0, $T$ being an
IET, $\mu_P(cod_P(\Discaccr))=0$, so $\mu_P(\sigma^{-1}(cod_P(\Discaccr)))=0$
too.
Hence $\mu_P(\sigma^{-1}(SP_{L_P})) = 0$, and again by measure preservation,
$\mu_P(SP_{L_P}) = 0$.  
\EPR

Then, given a monotonicity partition $P$, Lemma~\ref{fundlemma} applies
to $\Lang_{P}$ measured by $\mu_P$, so that Theorem~\ref{main1} applies to
$\Lang_{P}$ too, and we can get a conjugate IET~$\TietP$.

\vspace*{-0.2cm}
\begin{Prop}\label{iet_gives_the_same_iet_prop}
  Let $T$ be a piecewise increasing and minimal IET.  Let $P$ be a monotonicity
  partition of $I$, and let $\Lang_{P}$ be the associated shift
  with its conjugate IET $\TietP$. Then $\TietP=T$, i.e.:
\vspace*{-0.3cm}
\[\begin{array}{ccccccc}
\mbox{For every monotonicity partition } P, &  
             T & \stackrel{cod_P}{\longrightarrow} &\Lang_P & 
                 \stackrel{\mbox{\scriptsize Theorem~1}}{\longrightarrow} 
               & \TietP = T.
\end{array}\]
\vspace*{-0.7cm}
\end{Prop}
\PR In Equality~(\ref{pushf_basic_eq1}), $T(x)$ can replace $x$, and
by using~(\ref{conjug_basic_eq1}) we get for all $x \in I$, 
\vspace*{-0.3cm}
$$T(x) = \phi_{\mu_P}(cod_P(T(x))) =   
  \phi_{\mu_P}(\sigma (cod_P(x))).$$
Also, right-continuity of $cod_P$ (see Lemma~\ref{cod_increasing_lem}) means
that for every $x \in I$, $cod_P(x)$ is the infimum of $\{cod_P(y) \;|\; y>x,
\;y\in I\}$. Thus $cod_P(x)$ cannot be the right endpoint of a cylinder in
$\Lang_{P}$, hence $cod_P(x) \in \phi^{-1}(\iota(I))$.
In fact, since $\phi_{\mu_P}=\kappa \circ \phi$, and since $\iota \circ \kappa$
is the identity on $\phi(cod_P(I))$,
we get $cod_P(x) = {\phi}^{-1}(\iota(x))$ by composing
Equality~(\ref{pushf_basic_eq1}) to the left with $\phi^{-1} \circ \iota$.
By Proposition~\ref{conjug_prop}(d), $\TietP(\phi_{\mu_P}(w)) =
\phi_{\mu_P}(\sigma(w))$ on $\phi^{-1}(\iota(I))$, thus
with $w=\phi^{-1}(\iota(x))$, for all $x \in I$, 

\vspace*{-0.3cm}
$$\TietP(x) = {\phi_{\mu_P}}(\sigma(cod_P(x))),$$
\vspace*{-0.5cm}

\noindent
whence $T$ and $\TietP$ coincide on~$I$.
\EPR

Of course the above result applies to minimal {\em finite} IETs, as these are
piecewise increasing by definition.  In this case, the partition $P$ can be the
one induced by the set~$\Disc$ of discontinuities, or by any finite $B \supseteq
\Disc$.  Then, $cod_{P}$ becomes the usual {\smbf natural coding} of the orbits
of~$T$ with $B$ as set of separation points~\cite{Kea75}, that is, a coding
based on the assignment of a distinct letter to each component of a partition
on which $T$ is continuous in the interior of these components:

\begin{Cor}\label{finite_iet_gives_the_same_finite_iet_cor}
  Let $T$ be a finite minimal IET. Let $\Lang$ be the associated shift given by
  a natural coding of $T$,
  with its conjugate IET $\Tiet$.  Then $T = \Tiet$.
\end{Cor}

\noindent
Note that in this finite IET case, $T$ is reobtained as~$\Tiet$ from $\Lang$,
but according to the construction behind Lemma~\ref{fundlemma}, the set of
separation points~$\Biet$ is infinite, whereas~$\Disc \subset \Biet$ is finite
(see Example~\ref{iet-graphs}~(iv) further).

A shift $\Lang$ can also be recovered from $\Tiet$.
Recall that $P_A$ is the monotonicity partition of $I$ for $\Tiet$ relying on
the letters of~$A$ (see Lemma~\ref{rightcontinuous_lem}):
\begin{Prop}\label{coding_incr_prop}
Let $\Lang$ be any shift with a conjugate IET $\Tiet$ (Theorem~\ref{main1}
applies to~$\Lang$).  
Then for every $x \in I$, $cod_{P_A}(x) = \phi^{-1}(\iota(x))$, so that
$\Lang_{P_A}=\Lang$, i.e.:
\vspace*{-0.4cm}
\[\begin{array}{ccccccc}
\mbox{For the monotonicity partition } P_A, &  
         \Lang & \stackrel{\mbox{\scriptsize Theorem~1}}{\longrightarrow} & \TietPA
               & \stackrel{cod_{P_A}}{\longrightarrow} &\Lang_{P_A}=L.
\end{array}\]
\vspace*{-0.8cm}
\end{Prop}
\PR For $x \in I$, by definition of $cod_{0,P_A}$
and $P_A$, 
$cod_{P_A}(x) = a$ when $x \in I_{L,a}$, since $\alpha(I_{L,a})=a$.
Next, by definition of $I_{L,a}$, $\phi^{-1}(\iota(x)) \in Cyl(a)$, that is, the
first letter of $cod_{P_A}(x)$ and of $\phi^{-1}(\iota(x))$ are equal.
To check equality on their second letters, that is, the first letter of
$\sigma(cod_{P_A}(x))$ and of $\sigma(\phi^{-1}(\iota(x)))$,
consider that since $x$ above can be arbitrary, 
the first letter of $cod_{P_A}(\Tiet(x))$ and of $\phi^{-1}(\iota(\Tiet(x)))$
are equal too.
Then by Equality~(\ref{conjug_basic_eq1}),
$cod_{P_A}(\Tiet(x))= \sigma(cod_{P_A}(x))$,
and by Proposition~\ref{conjug_prop}~(a)(b),
$\phi^{-1}(\iota(\Tiet(x)))= \sigma(\phi^{-1}(\iota(x)))$.
By inductively applying the same arguments, equality holds for the next
letters. Hence $cod_{P_A}(I) = \phi^{-1}(\iota(I))$.~\EPR

If an IET is given with an infinite set $B \supseteq \Disc$ of separation
points, its natural coding is such that 
the coding alphabet must be then infinite too. However, in the case of $\Tiet$
with $\Biet$ as set of separation points, this coding is related to $\Lang$,
that is, to a language over a finite alphabet.
Indeed, let~$A_\infty$ be an infinite alphabet, and let $\alpha_\infty$ be a map
assigning a distinct letter of~$A_\infty$ to each component of the
partition~$PART_{\Lang,I}$ (cf.~p.~\pageref{partI_def}), so that its natural
coding is defined as 
\[\begin{array}{llll} 
 cod_{0,\infty}:\; & I \;\;\;\; & \rightarrow \;\;\;\; & A_\infty \\
          & x          & \mapsto &\left\{
      \begin{array}{lll} 
          \alpha_\infty(\Ikiet) & \mbox{\em if $x \in \Ikiet$}\\
          \alpha_\infty(x) & \mbox{\em if $x \in B_{\SP}$.}
       \end{array} \right.
\end{array}\]
\[\begin{array}{lllllll} 
 cod_\infty:\; & I \;\;\;\; & \rightarrow \;\;\;\; & A_\infty^\mathbb{N} \\
         & x & \mapsto &  cod_{0,\infty}(x)\;\;cod_{0,\infty}(\Tiet(x))\;\;cod_{0,\infty}(\Tiet^2(x))...
\end{array}\]

\noindent
Recall now that the partition $PART_\Lang$ of $\Lang$
(cf.~p.~\pageref{partL_def}) is made of the cylinders determined by the finite
words $v^{(k)}$ in $Fact_\Lang$, and by the infinite words in $W_{\SP}$.
We consider then the following injective recoding map of the letters of
$A_\infty$:
\[\begin{array}{llll}\label{zeta_infty_def}
 \zeta_\infty:\; & A_\infty & \rightarrow & Fact_\Lang \cup W_{\SP} \\
                    & a & \mapsto & \left\{
      \begin{array}{lll}
           v^{(k)} & \mbox{\em if $a=\alpha_\infty(\Ikiet)$} \\ 
           w &       \mbox{\em if $a=\alpha_\infty(x)$, where  
                           $x \in B_{\SP}$ and $\phi^{-1}(\iota(x))=w$.} 
\end{array} \right.
\end{array}\]

\noindent
Next, we define 

\vspace*{-0.5cm}
\[\begin{array}{lllllll} 
 \delta_{0,\infty}:\; & \Lang \;\;\;\; & \rightarrow \;\;\;\; & A_\infty,\\
\label{delta_0_infty_def}
\end{array}\]
\vspace*{-0.8cm}

\noindent
where for each $w \in \Lang$:
\begin{itemize}\itemsep=0cm
  \item $\delta_{0,\infty}(w) = (\zeta_\infty)^{-1}(au)$, if $au$ occurs as a
    prefix of $w$, with $a \in A$ and $u \in Fact_\Lang$ is non-left special,
    that is, by the construction of $PART_\Lang$ in Lemma~\ref{fundlemma}, $au$
    is a factor $v^{(k)}$ such that $Cyl_\Lang(v^{(k)}) \in PART_\Lang$.

 \item $\delta_{0,\infty}(w) = (\zeta_\infty)^{-1}(w)$, if $w$ is of the form
   $w=aw'$, with $a \in A$ and $w'$ is an infinite left special word in~$\SP$,
   that is, $w \in W_{\SP}$.
\end{itemize} 

\noindent
The map $\delta_{0,\infty}$ is well-defined because 
first, the two above cases cover $\Lang$ and are independent, reflecting
how $PART_\Lang$ is built;
second, if $au$ exists, it is unique since the $v^{(k)}$'s determining the
cylinders of $PART_\Lang$ can never be prefixes to each other (these cylinders
do not overlap).
Then, $\delta_{0,\infty}$ is extended as follows:
\[\begin{array}{lllllll} 
 \delta_\infty:\; & \Lang \;\;\;\; & \rightarrow \;\;\;\; & A_\infty^\mathbb{N} \\
           & w  & \mapsto &  \delta_{0,\infty}(w)\;\;\delta_{0,\infty}(\sigma(w))\;\;\delta_{0,\infty}(\sigma^2(w))...
\end{array}\]

\noindent
We can now map the words in $\Lang$ to symbolic orbits of $\Tiet$ obtained by
$cod_\infty$:
\begin{Prop}\label{codings_equiv_prop}
  Let $\Lang$ be any shift with a conjugate IET $\Tiet$ (Theorem~\ref{main1}
  applies to~$\Lang$).  Then for every $x \in I$, $cod_\infty(x) =
  \delta_\infty(\phi^{-1}(\iota(x)))$, i.e.:
\vspace*{-0.1cm}
\[\begin{array}{ccccccc}
    \Lang  
    & \stackrel{\mbox{\scriptsize Theorem~1}}{\longrightarrow} & \Tiet
    & \stackrel{cod_\infty}{\longrightarrow} & 
             cod_\infty(I) = \delta_\infty(cod_{P_A}(I)).
\end{array}\]
\vspace*{-0.4cm}
\end{Prop}
\PR 
For $x\in I$, the first letter of $cod_{\infty}(x)$ and of $\delta_{\infty}(w)$
where $w=\phi^{-1}(\iota(x))$, are equal.
Indeed, since $w$ belongs to some $Cyl(v^{(k)})$ or to~$W_{\SP}$,
$\delta_{0,\infty}(w)$ recovers the letter assigned by $\alpha_\infty$ to the
component of $PART_{\Lang,I}$ containing~$x$, i.e.~$cod_{\infty}(x)$.
Next, $\sigma (cod_\infty(x)) = cod_\infty(\Tiet(x))$, and
$\sigma(\delta_\infty(\phi^{-1}(\iota(x))))$ $=
\delta_\infty(\sigma(\phi^{-1}(\iota(x)))) = \delta_\infty(\phi^{-1}
(\iota(\Tiet(x))))$,
so that the second letters are also equal. By inductively applying the same
arguments, we obtain $cod_\infty(x) = \delta_\infty(\phi^{-1}(\iota(x)))$. Hence
$cod_\infty(I) = \delta_\infty(\phi^{-1}(\iota(I)))$, and by
Proposition~\ref{coding_incr_prop},
$\delta_\infty(\phi^{-1}(\iota(I)))=\delta_\infty(cod_{P_A}(I))$.
\EPR

\noindent
This mapping of $\Lang$ to the natural coding of $\Tiet$ will be in use in the
next section.
\section{The Linear Complexity Shift Case}

We now focus on shifts $\Lang$ with {\smbf linear complexity},
i.e.~those such that~$p_\Lang(n)=O(n)$, and prove Theorem~\ref{main2}.
For instance, aperiodic minimal shifts associated with primitive substitutions
have linear complexity~\cite{Pan84,Que10}, e.g.~the shift~associated with the
Thue-Morse substitution.
Since these shifts form a specific family of shifts with zero topological
entropy, the idea behind Theorem~\ref{main2} is to investigate the properties of
their conjugate IETs as given by Theorem~\ref{main1}.

\subsection{Almost Finite Interval Exchange Transformations}

We begin by defining a specific kind of IETs:

\vspace*{-0.1cm}
\begin{Def}\label{almost-finite-def}
An IET $T$ with a set $\Disc$ of discontinuities is {\smbf almost finite} if:
\vspace*{-0.1cm}
\begin{enumerate}\itemsep=-0.1cm
\vspace*{-0.1cm}
\item 
$T$ is piecewise increasing on $I$.

\item 
  The set $\Discacc$ of accumulation points of $\Disc$ is finite.

\item 
  The points in $\Disc$ belong to a finite number of full orbits of
  $T$\footnote{A sufficient condition for minimality of a finite IET is that the
    full orbits of the discontinuity points are infinite and
    distinct~\cite{Kea75} (cf.~the {\em infinite distinct orbit condition}
    (IDOC)).  An almost finite IET with an infinite number of discontinuities is
    thus far from satisfying the IDOC.}.
\end{enumerate}
\end{Def}

\vspace*{-0.1cm}
\noindent
In the above definition --~like in Definitions~\ref{iet_def}
and~\ref{iiet_def}~-- $\Disc$ can be replaced by a null-measure set of
separation points $B$ such that $B \supset \Disc$. Note also that a finite IET
is almost finite.
Now, the first property of almost finiteness is always satisfied by a
conjugate IET~$\Tiet$ (see Lemma~\ref{rightcontinuous_lem}).
About the second property, when $\Lang$ has linear complexity, we first show
the following technical result:
\begin{Lem}\label{finite_inf_sp_lem}
Let $\Lang$ be a minimal shift with linear complexity.  Then $\SP$ is
finite.
\end{Lem}
\PR For each finite $r \leq |\SP|$, there is $n_0>0$ such that for all $n\geq
n_0$, the number of distinct length-$n$ prefixes of the words in $\SP$ is at
least~$r$. 
If~$sp_{\Lang,l}(n)$ denotes the number of left special factors in
$Fact_\Lang(n)$, then $sp_{\Lang,l}(n) \ge r$.  Since $\Lang$ is
minimal, $Fact_\Lang$ is prolongable, and thus for all $n\geq n_0$,
$p_\Lang(n+1) - p_\Lang(n) \ge sp_{\Lang,l}(n)$.
Therefore, if $\SP$ was infinite, for every $r \in \mathbb{N}$, there would
exist $n_0>0$ such that for all $n\geq n_0$, $p_\Lang(n+1) - p_\Lang(n) \ge r$,
whence $p_L(n) \neq O(n)$.  \EPR

\noindent
Then: 
\begin{Lem}\label{finite_accum_lem} 
  Let $\Lang$ be an aperiodic minimal shift with linear complexity, and
  let~$\Tiet$ be its conjugate IET.  Then the set $\Discietacc$ of the
  accumulation points of the set $\Disciet$ of the discontinuities of~$\Tiet$ is
  finite.
\end{Lem}
\PR By definition, $\Disciet \subseteq \Biet$. Next, 
$\Bietacc$ is determined by $W_{\SP}$ since $\phi_\mu(W_{\SP}) = \Bietacc$,  
and we also have $W_{\SP}=\sigma^{-1}(\SP)$.
The alphabet $A$ is finite, and by Lemma~\ref{finite_inf_sp_lem} $\SP$ is too,
whence the result.
\EPR

In order to prove the third property of almost finiteness for $\Tiet$, the idea
is to study similar relationships to known ones for finite IETs, between
natural~codings, complexity of the associated shifts, and distinctiveness of the
orbits of the separation points~\cite{Kea75}.
However, as seen in Section~\ref{coding_reconstruct_sec}, a natural
coding~$cod_\infty$ of an infinite IET $T$ uses an infinite alphabet~$A_\infty$.
To deal with this situation, we shall code $T$'s orbits using finer and finer
{\em finite} partitions of~$I$, together with finite alphabets with more
and more letters,
as is developed in the next section.

\subsection{A Finite Number of Full Orbits}

For an IET~$T$ with $B \supseteq \Disc$ as set of separation points, let $\Bacc$
denote the set of accumulation points of $B$ in~$\overline{I}$,  
and let $\Bacc = \Baccone \sqcup \Bacctwo$, where $\Baccone$ denote the
one-sided accumulation points, and $\Bacctwo$ the two-sided ones.
An ordered finite subset $B_m=\{0 = b_0, b_{1}, \dots, b_{m-1}\}$, $m> 0$, of
$B$~is
called {\smbf admissible} if for $B_m \cup \{b_m=1\}$, each induced interval
$I_{m,i}=[b_{i}, b_{i+1})$ has one of the following~types:

\begin{enumerate}\itemsep=0cm
 \item $I_{m,i}$ is an interval of $T$ with respect to $B$ ($b_i$ is the only
   point of $B$ in $I_{m,i}$). 
 \item $I_{m,i}$ covers infinitely many consecutive intervals of $T$, 
   and either:
\vspace*{-0.2cm}
{\small
   \begin{enumerate}[(a)]\itemsep=0cm
      \item $[b_{i}, b_{i+1}]$ has exactly one point of $\Bacc$,
        belonging to~$\Baccone$, either $b_{i}$~or~$b_{i+1}$;
      \item $(b_{i}, b_{i+1})$ has exactly one point of $\Bacc$, 
        belonging to~$\Bacctwo$, with $b_{i},b_{i+1}\notin\Bacc$.
   \end{enumerate}
}
\end{enumerate}

\noindent
Thus, an admissible $B_m$ induces a finite partition $\bigsqcup_i I_{m,i}$ of
$I$, where for each point $x \in \Bietacc$, there is exactly one $I_{m,i_x}$
such that almost every point of any sequence of points of $B$ converging to $x$
lies in $I_{m,i_x}$.
For instance, here~is~a possible subset~$B_4$ of a set $B$ with two accumulation
points $x_1 \in \Baccone$ and $x_2 \in \Bacctwo$, such that the induced
intervals $I_{4,0},I_{4,2}$ are of type~1, $I_{4,1}$ of type~2a, $I_{4,3}$ of
type~2b:

\vspace*{-0.1cm}
\begin{center} 
\includegraphics[width=9.5cm]{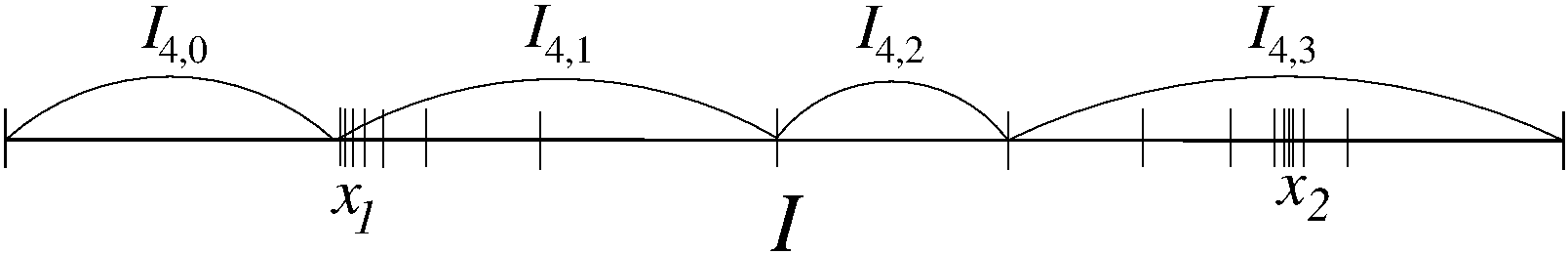}
\end{center}
\vspace*{-0.7cm}

\begin{Lem}\label{b_m_exist_lem} 
  Let $B$ be a set of separation points.
Then $\Bacc$ is finite iff for some $m>0$, there is an admissible~$B_m\subset
B$.
\end{Lem}
\PR 
$(\Leftarrow)$: Trivial. 
$(\Rightarrow)$: First, for each point $x \in \Bacctwo$ there is an interval of
type~2b containing it, because such an~$x$ is the intersection of a sequence of
closed nested intervals with endpoints in~$B$, and if every interval in the
sequence contained a point in $\Bacc$, then $\Bacc$ would be infinite.
Next, by the same argument, each $x \in \Baccone$ which is the intersection of a
sequence of nested closed intervals to its right is the left endpoint of an
interval of type~2a, while in its other side, $x$ is the right endpoint of an
interval of type~1, or else $x=0$. The same occurs the other way around when the
accumulation is on the left.
Then there are finitely many points of $B$ in the complement of the union of the
above intervals, and together with the endpoints of those intervals, they form
an admissible~$B_m$.\EPR

Given an admissible set $B_m$, let~$A_m$ be an alphabet with $m$ letters,
and let $\alpha_m$ be a map assigning a distinct letter in $A_m$ to each
$I_{m,i}$ induced by~$B_m$, so that
\vspace*{-0.1cm}
\[\begin{array}{llll} 
 cod_{0,m}:\; & I \;\;\;\; & \rightarrow \;\;\;\; & A_m \\
          & x    & \mapsto &\mbox{\em $\alpha_m(I_{m,i})$, \;\;\;if $x \in I_{m,i}$.} 
\end{array}\]
\vspace*{-0.2cm}
\[\begin{array}{lllllll} 
 cod_m:\; & I \;\;\;\; & \rightarrow \;\;\;\; & A_m^\mathbb{N} \\
   & x  & \mapsto &  cod_{0,m}(x)\;\;cod_{0,m}(T(x))\;\;cod_{0,m}(T^2(x))...
\end{array}\]

\noindent
The {\smbf associated shift} $\Lang_m$ of $cod_m$ is then as usual defined as
$\overline{cod_m(I)}$.
Now, let $\Lang$ be a shift over $A$, with a conjugate IET $\Tiet$ having
$\Biet$ as set of separation points, and $\Bietacc = \Bietaccone \sqcup
\Bietacctwo$ as accumulation point sets. Let~$\Bietm$ be an admissible subset of
$\Biet$, inducing a finite partition~$\bigsqcup_{i} \Iietmi$ of~$I$.
Then, 
we define a recoding map similar to and based on $\zeta_\infty$ (see
p.~\pageref{zeta_infty_def}), sending each letter of~$A_m$ to words over $A$,
which are either $v^{(k)}$'s or common prefixes of sets of $v^{(k)}$'s (the
$v^{(k)}$'s determining the cylinder components of~$PART_\Lang$):
\vspace*{-0.2cm}
\[\begin{array}{llll}
 \zeta_m:\; & A_m  \;\; & \rightarrow \;\;
                         Fact_\Lang \cup (Fact_\Lang \times Fact_\Lang), 
\end{array}\]
\vspace*{-0.7cm}

\noindent
where for each $a \in A_m$:
\vspace*{-0.2cm}
\begin{enumerate}\itemsep=0cm
 \item $\zeta_m(a) = \zeta_\infty(\alpha_\infty(\Iietmi))$, if
   $a=\alpha_m(\Iietmi)$ where $\Iietmi$ is of type~1,
   that is, $\zeta_m$ yields the same $v^{(k)}$ as $\zeta_\infty$ if $a$ is
   assigned to a type~1 interval. 
 \item \begin{enumerate}[(a)]\itemsep=0cm
   \item $\zeta_m(a)=\vmaxpref$, if $a=\alpha_m(\Iietmi)$ where $\Iietmi$ is of
     type~2a, and $\vmaxpref$ is~the longest common prefix of
     $\{\zeta_\infty(\alpha_\infty(J))| J\in PART_{\Lang,I}, J\subset\Iietmi\}$,
     that is, $\zeta_m$ yields the longest common prefix of the $v^{(k)}$'s
     given by $\zeta_\infty$ applied to all the intervals covered by $\Iietmi$.

  \item $\zeta_m(a)=(\vmaxprefl, \vmaxprefr)$, if $a=\alpha_m(\Iietmi)$ where
    $\Iietmi$ of type~2b, and $\vmaxprefl, \vmaxprefr$ are defined as follows:
    let $\Iietmi= [b_i,b_{i+1})$ be divided into $J_l=[b_i,x)$ and
        $J_r=[x,b_{i+1})$, where $x \in \Bietacctwo$ is the accumulation point
          in~$\Iietmi$; then $\vmaxprefl$ (resp. $\vmaxprefr$) is the longest
          common prefix of $\{\zeta_\infty(\alpha_\infty(J))\;|\; J \in
          PART_{\Lang,I}, J\in J_l\}$ (resp. $J \in J_r$).
    \end{enumerate}
\end{enumerate}

\noindent
The case 2(b) above deals with the possible existence of points~$x \in
\Bietacctwo$ such that $\phi_\mu^{-1}(x)$ is made of two consecutive distinct
words in $W_{\SP}$.

We show then that if $m$ is sufficiently large, the map $\zeta_m$
can be used like~$\zeta_\infty$ to define a map $\delta_m$ similar to
$\delta_\infty$.
First of all, we say a sequence $\{B_m\}_{m \geq m_0}$, with $m_0>0$, is
{\smbf admissible} if for all $m \geq m_0$: (i) $B_m$ is an admissible subset
of~$B$; (ii)~$|B_{m+1}| = |B_m| + 1$; (iii)~$B = \bigcup_{m \geq m_0} B_m \cup
(\Bacc \setminus \{1\})$.

\begin{Lem}\label{admiss_seq_lem}
Let $B$ be a set of separation points. 
Then $\Bacc$ is finite iff for some $m_0>0$, there is an admissible
sequence $\{B_m\}_{m \geq m_0}$.
\end{Lem}
\PR 
$(\Leftarrow)$: A consequence of Lemma~\ref{b_m_exist_lem}.
~$(\Rightarrow)$: By Lemma~\ref{b_m_exist_lem} also, for some $m_0>0$ there~is an
admissible $B_{m_0}$ of $B$.
If the first point $x$ of $\Bacc$ in $I$ lies in an interval $I_{m_0,i}$ of
type~2b, we add to $B_{m_0}$ the smallest point of~$B$ in the interior of
$I_{m_0,i}$ to obtain $B_{m_0+1}$, then the greatest one to obtain
$B_{m_0+2}$.
If~$x$ lies in $I_{m_0,i}$ of type~2a we do the same, but using the smallest or
the greatest point of $B$ in the interior of $I_{m_0,i}$, according to if the
point is an accumulation point on the left or the right.  Then we proceed to the
next point of $\Bacc$, and so on, going back to the first one after having
processed the last one.  
\EPR

\noindent
We denote by~$F_{m,1}$ (resp. $F_{m,2}$) the set of the $v^{(k)}$'s (resp. the
prefixes of the~$v^{(k)}$'s) given by $\zeta_m(\alpha_m(\Iietmi))$ when
$\Iietmi$ is of type~1 (resp. type~2, including $\vmaxprefl$ and $\vmaxprefr$
as distinct words when $\Iietmi$ is of type~2b).
The technicalities in the proof of the following lemma come from the topological
differences between $\Lang$ and~$I$:
\begin{Lem}\label{zeta_m_inject}
Let $\Lang$ be a shift with a conjugate $\Tiet$ and
an admissible sequence $\{\Bietm\}_{m \geq m_0}$, $m_0>0$.  Then there is
$m_1\geq m_0$, such that for all $m \geq m_1$, the map $\zeta_m$ is injective,
such that for $i=1,2$, $F_{m,i}$ does not include the empty word, nor words
prefix to each other, and no word in $F_{m, 1}$ is a prefix of a word in $F_{m,
  2}$.
\end{Lem}
\PR By the construction of $PART_\Lang$ in Lemma~\ref{fundlemma}, the set
$F_{m,1}$, being made of $v^{(k)}$'s, has the first two required properties. 
Next, a word~$v$ in $F_{m,1}$ cannot be a prefix of any word in~$F_{m,2}$, since
$v$ would be then a common prefix of words of type~$v^{(k)}$, and again the
$v^{(k)}$'s are not prefixes to each other.

Now, from the point-of-view of $\Lang$, each word $w \in W_{\SP}$ is an
accumulation of endpoints of cylinder components of $PART_\Lang$, and can be
seen as the intersection of a sequence of nested cylinders of~$\Lang$
(see again the proof of Lemma~\ref{fundlemma}).
Also, by Lemma~\ref{admiss_seq_lem},
the set~$\Bietacc$ must be finite, so that $W_{\SP}$ is too.
Therefore, there exists~$n>0$ such that at the $n$th step of the construction of
$PART_\Lang$, the nested cylinders to be refined at each further step are such
that they contain only one word in~$W_{\SP}$, and such that they are all
separated by non-empty cylinders. Thus for each $w \in W_{\SP}$, there is a
sequence of nested cylinders $\{C_{w,j}\}_{j \in \mathbb{N}^*}$ such that
$\bigcap_{j\in \mathbb{N}^*}C_{w,j} = \{w\}$, separated 
from the other sequences by non-empty cylinders of $\Lang$,  
and having one of these two forms:

\vspace*{-0.2cm}
\begin{enumerate}[a$_\Lang$)]\itemsep=-0.1cm
\item $\{C_{w,j}\}_{j \in \mathbb{N}^*}$ is such that $w$ is the smallest or the
  greatest element of every~$C_{w,j}$,
so that $w$ is a one-sided accumulation word in $\Lang$.

\item $\{C_{w,j}\}_{j \in \mathbb{N}^*}$ is such that $w$ belongs to the
  interior of every cylinder~$C_{w,j}$,
  so that $w$ is a two-sided accumulation word in $\Lang$.
\end{enumerate}
\vspace*{-0.2cm}

From the point-of-view of~$I$, the corresponding sequences
$\{\phi_\mu(C_{w,j})\}_{j \in \mathbb{N}^*}$~are pairwise disjoint sequences of
nested closed intervals in $I$,
such that $\bigcap_{j \in \mathbb{N}^*}\phi_\mu(C_{w,j})$ $=\{\phi_\mu(w)\}$,
with $\phi_\mu(w) \in \Bietacc$,
taking one of these two forms:

\vspace*{-0.2cm}
\begin{enumerate}[a$_{I}$)]\itemsep=-0.1cm
\item $\{\phi_\mu(C_{w,j})\}_{j \in \mathbb{N}^*}$ is such that 
  $\phi_\mu(w)$ is equal to the right or the left endpoint of every interval
  $\phi_\mu(C_{w,j})$ so that $\phi_\mu(w)\in \Bietaccone$ or $\Bietacctwo$.
\item $\{\phi_\mu(C_{w,j})\}_{j \in \mathbb{N}^*}$ is such that
  $\phi_\mu(w)$ belongs to the interior of $\phi_\mu(C_{w,j})$, so that
  $\phi_\mu(w) \in \Bietacctwo$.
\end{enumerate}
\vspace*{-0.2cm}

From the point-of-view of $\{\Bietm\}_{m \geq m_0}$, for each $m \geq m_0$, every
$x \in \Bietacc$ is such that there is a~$\Iietmi$ of type~2 containing almost
every point of any sequence in~$\Biet$ converging to $x$,
so that putting $J_{x,m} = \overline{\Iietmi}$, we get a sequence of nested
intervals $\{J_{x,m}\}_{m \geq m_0}$, with $\bigcap_{m \geq m_0}J_{x,m}=\{x\}$,
taking one of these three forms:

\vspace*{-0.2cm}
\begin{enumerate}[a$_{B}$)]\itemsep=-0.1cm
\item $\{J_{x,m}\}_{m \geq m_0}$ is made of type~2a intervals, i.e.~$x \in
  \Bietaccone$, where $\phi_\mu^{-1}(x)=w$, $w \in W_{\SP}$. There is then a
  corresponding sequence $\{C_{w,j}\}_{j \in \mathbb{N}^*}$ of form~a$_\Lang$
  with intersection $\{w\}$, and $\{\phi_\mu(C_{w,j})\}_{j \in \mathbb{N}^*}$ of
  form~a$_I$ with intersection $\{x\}$.

 \item $\{J_{x,m}\}_{m \geq m_0}$ is made of type~2b intervals, i.e.~$x \in
   \Bietacctwo$, where either: 
\vspace*{-0.2cm}
  \begin{enumerate}[i)]\itemsep=-0.1cm
    \item $\phi_\mu^{-1}(x)=w$, with $w \in W_{\SP}$. There is then a
      corresponding sequence $\{C_{w,j}\}_{j \in \mathbb{N}^*}$ of
      form~b$_\Lang$ with intersection $\{w\}$, and $\{\phi_\mu(C_{w,j})\}_{j
        \in \mathbb{N}^*}$ of form~b$_I$ with intersection $\{x\}$.

    \item $\phi_\mu^{-1}(x)=\{w,w'\}$, with $w,w' \in W_{\SP}$. There are then
      corresponding sequences $\{C_{w,j}\}_{j \in \mathbb{N}^*}$,
      $\{C_{w',j}\}_{j \in \mathbb{N}^*}$ of form~a$_\Lang$, respectively with
      intersection $\{w\}$ and $\{w'\}$, and $\{\phi_\mu(C_{w,j})\}_{j \in
        \mathbb{N}^*}$, $\{\phi_\mu(C_{w',j})\}_{j \in \mathbb{N}^*}$ of
      form~a$_I$, both with intersection $\{x\}$.
  \end{enumerate}
\vspace*{-0.1cm}
\end{enumerate}
\vspace*{-0.2cm}

\vbox{
Now, since the sequences $\{\phi_\mu(C_{w,j})\}_{j \in \mathbb{N}^*}$ and
$\{J_{x,m}\}_{m \geq m_0}$ have the same intersection point in $\Bietacc$, and
because of the properties of the $C_{w,1}$'s, there exists $m_1>0$, such that
one of the following cases occur:

\vspace*{-0.2cm}
\begin{itemize}\itemsep=-0.1cm
\item $\{J_{x,m}\}_{m \geq m_0}$ is of form~a$_B$ or~b$_B$i, and for all
  $m>m_1$, $J_{x,m} \subset \phi_\mu(C_{w,1})$, that is, the word determining
  the cylinder $C_{w,1}$ is a common prefix of the words in
  $\{\zeta_\infty(\alpha_\infty(J_{x,m}))\mid m>m_1\}$, i.e.~a prefix of the
  corresponding $\vmaxprefgen \in F_{m,2}$.

\item $\{J_{x,m}\}_{m \geq m_0}$ is of form~b$_B$ii, and for all $m>m_1$,
  $J_{x,m} \subset \phi_\mu(C_{w,1})\cup \phi_\mu(C_{w',1})$,
  that is, the two words determining $C_{w,1}$ and $C_{w',1}$ are respective
  common prefixes of the factors in $\{\zeta_\infty(\alpha_\infty(J_l))\mid
  m>m_1\}$ and $\zeta_\infty(\alpha_\infty(J_r))\mid m>m_1\}$, where $J_l$
  and~$J_r$ are the two components of $J_{x,m} \setminus x$, i.e.~prefixes of
  the corresponding $(\vmaxprefl, \vmaxprefr)$, with $\vmaxprefl, \vmaxprefr \in
  F_{m,2}$.
\end{itemize}
\vspace*{-0.2cm}

\noindent
Thus, for all $m \ge m_1$, since all the $C_{w,1}$'s are distinct, the set
$F_{m,2}$ does not contain the empty word, nor words which are prefix to each
other. 
\EPR
}

Similarly to $\delta_{0,\infty}$, we can now define $\delta_{0,m}$ for each
$m>m_1$, using the properties of~$\zeta_m$, so as to obtain letters of~$A_m$
out of the prefixes of the words in~$\Lang$:
\[\begin{array}{llll} 
 \delta_{0,m}:\; & \Lang \;\;\;\; & \rightarrow \;\;\;\; & A_m, \\
\end{array}\]
\vspace*{-0.6cm}

\noindent
where for each $w \in \Lang$: 
\begin{enumerate}\itemsep=-0.1cm
\item $\delta_{0,m}(w)=(\zeta_m)^{-1}(au)$, if $au$ occurs as a prefix of $w$ in
  $F_{m,1}$, that is, $au$ is a factor $v^{(k)}$ which belongs to
  $\zeta_m(A_m)$, and such that $Cyl_\Lang(v^{(k)}) \in PART_\Lang$.
\item Or, one of the two following cases occurs:
\vspace*{-0.2cm}
  \begin{itemize}\itemsep=0cm
   \item $\delta_{0,m}(w)=(\zeta_m)^{-1}(\vmaxprefgen)$, if $\vmaxprefgen$
     occurs as a prefix of $w$ in $F_{m,2}$, which belongs as a single word to
     $\zeta_m(A_m)$.
  \item $\delta_{0,m}(w)= (\zeta_m)^{-1}(\vmaxprefl, \vmaxprefr)$, if
    $\vmaxprefl$ or $\vmaxprefr$ occurs as a prefix of~$w$ in $F_{m,2}$, which
    belongs to a pair $(\vmaxprefl,\vmaxprefr)$ in $\zeta_m(A_m)$.
  \end{itemize}
\end{enumerate}

\noindent
This map is well-defined because 
first, all the above cases cover $\Lang$ by definition of $F_{m,1}$ and
$F_{m,2}$; 
second, the properties of $F_{m,1}$ and $F_{m,2}$ proved by
Lemma~\ref{zeta_m_inject} imply that taken in order, only one case applies for
each $w \in \Lang$.
Accordingly, $\delta_{0,m}$ is extended as follows to obtain words in
$A^{\mathbb{N}}_m$ from words in $\Lang$:
\[\begin{array}{lllllll} 
 \delta_m:\; & \Lang \;\;\;\; & \rightarrow \;\;\;\; & A_m^\mathbb{N} \\
             & w & \mapsto &  \delta_{0,m}(w)\;\;\delta_{0,m}(\sigma(w))\;\;\delta_{0,m}(\sigma^2(w))...
\end{array}\]
\vspace*{-0.3cm}

\noindent
We then have a similar result to Proposition~\ref{codings_equiv_prop} 
($m_1$ comes from~Lemma~\ref{zeta_m_inject}):

\begin{Lem}\label{codings_equiv_m_lem}
Let $\Lang$ be any shift with a conjugate $\Tiet$ and an admissible
sequence $\{\Bietm\}_{m \geq m_1}$.  Then, for all $m\geq m_1$ and for all $x
\in I$, $cod_m(x) = \delta_m(\phi^{-1}(\iota(x)))$, i.e.:
\vspace*{-0.2cm}
\[\begin{array}{ccccccc}
    \Lang  
    & \stackrel{\mbox{\scriptsize Theorem~1}}{\longrightarrow} & \Tiet
    & \stackrel{cod_m}{\longrightarrow} & 
             L_m = \delta_m(L).
\end{array}\]
\vspace*{-0.5cm}
\end{Lem}
\PR 
For $x\in I$, the first letter of $cod_m(x)$ and of $\delta_m(w)$ where
$w=\phi^{-1}(\iota(x))$, are equal.
Indeed, 
$w$ belongs to some cylinder determined by a word in $F_{m,1} \cup
F_{m,2}$. Then, Case~(1) of $\delta_{0,m}(w)$'s definition applies iff~$x$
belongs to a type~1 interval, and $\delta_{0,m}(w)$ recovers the letter assigned
by $\alpha_m$ to this interval, i.e.~$cod_m(x)$. Otherwise, $x$ belongs to a
type~2 interval, and Case~(2) of $\delta_{0,m}(w)$'s definition necessarily
applies, with the same effect.
This is well-defined since by Lemma~\ref{zeta_m_inject} only words in $F_{m, 2}$
can be prefixes of words in $F_{m, 1}$.
Next, similarly to Proposition~\ref{codings_equiv_prop}, using the conjugacies
given by Proposition~\ref{conjug_prop}, the same is true for the other letters
of $cod_m(x)$ and $\delta_m(w)$. Hence $cod_m(I) =
\delta_m(\phi^{-1}(\iota(I)))$.~\EPR

\begin{Lem}\label{samecomplexity_lem}
Let $\Lang$ be any shift with a conjugate $\Tiet$ and an
admissible sequence $\{\Bietm\}_{m \geq m_1}$. 
Then, for all $m \geq m_1$, there exists $h_m>0$ such that for all $n > 0$,
$p_{\Lang_m}(n) \leq p_{\Lang}(n+h_m)$,
\end{Lem}
\PR Fix some $m\geq m_1$, and let $w \in \Lang$ and $w'=\delta_m(w) \in
\Lang_m$. Let $h_m$ be the maximum length of the factors in $F_{m,1} \cup
F_{m,2}$.  According to the definition of~$\delta_{0,m}$ and to
Lemma~\ref{codings_equiv_m_lem}, for each $n>0$, the length-$(n+h_m)$ prefix of
$w$ determines by~$\delta_m$ at least one length-$n$ prefix of $w'$,
hence, $p_{w'}(n) \leq p_{w}(n+h_m)$.
Now, $\Lang$~being minimal, $\Lang_m$ is too, 
and $w,w'$ are minimal words.  Thus each of the factors in $Fact_w$,
i.e.~$Fact_{\Lang}$, occurs as a prefix in the words in $\{\sigma^i(w)\}_{i\in
  \mathbb{N}}$, and these factors determine $Fact_{w'}$,
i.e.~$Fact_{\Lang_m}$. Hence for all $n>0$, $p_{\Lang_m}(n) \leq
p_{\Lang}(n+h_m)$.\EPR

\begin{Lem}\label{gogo_girl_lem} 
Let $\Lang$ be an aperiodic minimal shift with linear complexity.  Let $\Tiet$
be its conjugate IET.  Then~$\Biet$ belongs to a finite number of full orbits of
$\Tiet$.
\end{Lem}
\PR According to the proof of Lemma~\ref{finite_accum_lem}, $\Bietacc$ is
finite, thus according to Lemma~\ref{admiss_seq_lem} there is an admissible
sequence $\{\Bietm\}_{m \geq m_1}$ of $\Biet$, with~$m_1>0$ given by
Lemma~\ref{zeta_m_inject}. Fix some $m \geq m_1$, and consider the shift
$\Lang_m$ associated with $cod_m$.
For each $v \in Fact_{\Lang_m}$, let $E_v = \{x \in I\;|\; cod_{0,m}(x)
cod_{0,m}(\Tiet(x)) \dots$ $cod_{0,m}(\Tiet^{|v|-1}(x)) = v\}$,
i.e.~the set of the points in $I$ having symbolic orbits by $cod_m$ with~$v$ as
prefix
(note that $E_v$ is not necessarily connected).  For each $n>0$, $\bigsqcup_{v
  \in Fact_{\Lang_m}(n)} E_v$ is a partition of~$I$, and the set
$\{\Tiet^{-j}(\Bietm)\;|\; j=0,...,n-1\}$ defines the endpoints of the interval
components of this partition.
Then, for any point $x \in \Bietm$, if the set $\Tiet^{-n}(x)$
has a point $y$ in the interior of some~$E_v$ with~$v \in Fact_{\Lang_m}(n)$, it
means that to each side of~$y$ there are components of $E_{v_1}$ and $E_{v_2}$,
with $v_1, v_2 \in Fact_{\Lang_m}(n+1)$ with $v$ as prefix, and $v_1 \neq v_2$.
Thus~$v$ is a right special factor, and it induces an increase in the
complexity~$p_{\Lang_m}$ of $\Lang_m$.

Now, let $k_m \leq |\Bietm|=m$ be the number of distinct full orbits of
$\Bietm$. 
Since $\Lang_m$ is aperiodic minimal,
each of these orbits may contain at most one occurrence of each point in
$\Bietm$.
And since $\Bietm$ is finite, in each of these orbits there is a least $-j \in
\mathbb{Z}^-$ for which a point of $\Bietm$ occurs in $\{\Tiet^{-j}(\Bietm) \mid
j=0,...\}$.
Thus there are at least~$k_m$ points in $\Bietm$
such that for all $n>0$, their preimages by~$\Tiet^{-n}$ avoid $\Bietm$, and
being in the interior of some $E_v$, they induce an increase of the number of
the length-$n$ factors.
Hence, for all $m \geq m_1$ and for all $n>0$, $k_m n \leq p_{\Lang_m}(n)$. 
By Lemma~\ref{samecomplexity_lem}, and since $\Lang$ has linear
complexity, there exists $a>1$ such that 
\vspace*{-0.2cm}
$$ p_{\Lang_m}(n) \leq p_{\Lang}(n + h_m) \leq a(n + h_m).$$
that is, 
\vspace*{-0.2cm}
$$ k_m \leq a + \frac{ah_m}{n}.$$
Since $n$ is a free variable, we have $k_m \leq a$ for all $m \geq m_1$, and
thus letting $m$ go to infinity, this inequality holds also for $\Biet$
of~$\Tiet$.
\EPR

We can then prove Theorem~\ref{main2} as stated in the introduction:

\PRnotheo Since $\Lang$ has linear complexity, Theorem~\ref{main1} applies, and
we get a conjugate IET $\Tiet$.  By Lemmas~\ref{rightcontinuous_lem}
and~\ref{finite_accum_lem}, $\Tiet$ has the two first properties of almost
finiteness. Next, applying Lemma~\ref{gogo_girl_lem}, and using the fact that
$\Disciet \subseteq \Biet$, then~$\Tiet$ has also the third property of almost
finiteness.  \EPR

\section{Building Infinite Interval Exchange Transformations}

In this section, we present results in the context of zero entropy about how to
build effective approximations of the conjugate IETs $\Tiet$ given by
Theorem~\ref{main1} for a set of corresponding $\sigma$-invariant measures on
the shifts $\Lang$.

\subsection{An Approximation Scheme}
Given a shift $\Lang$ over an ordered alphabet $A$, and its complexity $p_\Lang$,
we define a sequence of maps $\{T_n\}_{n>1}$, 
where each $T_n$ is a map $I \rightarrow I$, $x \mapsto y=T_n(x)$, such that
the source on the $x$-axis is divided into $p_{\Lang}(n)$ right-open intervals
of equal length, and the range on the $y$-axis is divided into $p_{\Lang}(n-1)$
ones. The intervals of the source are then put into correspondence with the
cylinders of the factors in $Fact_{\Lang}(n)$, using the lexicographic order
between these factors, and the same is done for the intervals of the range with
cylinders of the factors in $Fact_{\Lang}(n-1)$.
Next, $T_n$ is defined as the piecewise affine map which sends for each $v \in
Fact_\Lang(n)$ the interval corresponding to $Cyl_\Lang(v)$ to the interval
corresponding to $Cyl_\Lang(\sigma(v))$, where $\sigma(v)$ denotes~$v$
minus its first letter,
by using a slope $\frac{p_{\Lang}(n)}{p_{\Lang}(n-1)}$.

Here is first a technical result about the above slopes converging
to~$1$. 
Recall that~$sp_{\Lang,l}(n)$ denotes the number of left special factors
in $Fact_{\Lang}(n)$, and let~$sp_{\Lang,r}(n)$ denote the same for the right
special~factors:
\begin{Lem}\label{slopes_lem}
   Let $\Lang$ be a shift such that $Fact_\Lang$ is prolongable. 
   Let $\{n_s\}_{s \in \mathbb{N}}$ be a subsequence of
   $\mathbb{N}\setminus\{0,1\}$. Then $\lim_{s\rightarrow \infty}
        \frac{p_{\Lang}(n_s)}{p_{\Lang}(n_s-1)} = 1$ iff $\lim_{s\rightarrow
          \infty}\frac{sp_{\Lang,l}(n_s-1)}{p_{\Lang}(n_s-1)} = 0$.
  The same holds for $sp_{\Lang,r}(n_s-1)$.
\end{Lem}
\PR We prove the result for $sp_{\Lang,l}$ (the proof for $sp_{\Lang,r}$ is
similar):\\ 
($\Rightarrow$): 
$\frac{p_{\Lang}(n_s)}{p_{\Lang}(n_s-1)} - 1 =
\frac{p_{\Lang}(n_s) - p_{\Lang}(n_s-1)}{p_{\Lang}(n_s-1)} \ge
\frac{sp_{\Lang,l}(n_s-1)}{p_{\Lang}(n_s-1)} \ge 0$.

\noindent
($\Leftarrow$):
$\frac{|A|\;sp_{\Lang,l}(n_s-1)}{p_{\Lang}(n_s-1)} \ge
\frac{p_{\Lang}(n_s)-p_{\Lang}(n_s-1)}{p_{\Lang}(n_s-1)} =
\frac{p_{\Lang}(n_s)}{p_{\Lang}(n_s-1)} -1 \ge 0$
\EPR

\noindent
For zero entropy,
using subsequences is indeed useful:

\begin{Lem}\label{slopesentropy0_lem}
   Let $\Lang$ be a minimal aperiodic shift with zero topological entropy. Then
   there is a subsequence $\{n_s\}_{s \in \mathbb{N}}$ of
   $\mathbb{N}\setminus\{0,1\}$ such that $\lim_{s \rightarrow \infty}
   \frac{p_{\Lang}(n_s)}{p_{\Lang}(n_s-1)} = 1$.
\end{Lem}
\PR Since $\Lang$ is minimal aperiodic, then for all $n>1$,
$\frac{p_{\Lang}(n)}{p_{\Lang}(n-1)}>
1$. Next,~$\frac{p_{\Lang}(n)}{p_{\Lang}(n-1)}$ is not bounded away from~1.
Assume~on~the contrary there is $c > 0$ such that for all $n>1$, $p_{\Lang}(n)
\ge (1+c)p_{\Lang}(n-1)$.
Hence $p_{\Lang}(n) \ge (1+c)^{n}$, and the topological entropy of $\Lang$
would not be less than $\log(1+c)$, a contradiction.\EPR

Let $Fact_{\Lang}(u,n)=\{v \in Fact_{\Lang}(n) \;|\; v = uv', v' \in A^*\}$,
i.e.~the set of~length-$n$ factors with $u$ as prefix, and let
$p_\Lang(u,n)=|Fact_{\Lang}(u,n)|$.  In the next result we simultaneously
build measures on $\Lang$ and approximations of the corresponding~$\Tiet$:

\begin{Prop}\label{approx_prop}
Let $\Lang$ be a minimal aperiodic shift with zero topological entropy. Let
$\{T_n\}_{n>1}$ be the sequence of maps as defined above, and let
$\Tseq=\{T_{n_s}\}_{s \in \mathbb{N}}$, where
$\{n_s\}_{s \in \mathbb{N}}$ is a subsequence of $\mathbb{N}\setminus\{0,1\}$
such that $\lim_{s \rightarrow \infty} \frac{p_{\Lang}(n_s)}{p_{\Lang}(n_s-1)} =
1$. Then:
\vspace*{-0.2cm}
\begin{enumerate}\itemsep=-0.1cm
   \item There is at least one shift-invariant measure $\mu$ on $\Lang$, induced
     by the $T_{n_s}$'s interval lengths.
   \item There is a subsequence of $\Tseq$ converging to a map~$T$, equal almost
     everywhere to the conjugate IET $\Tiet$ obtained by Theorem~\ref{main1}
     from $\Lang$ measured by~$\mu$.
\end{enumerate}
\vspace*{-0.3cm}
\end{Prop}
\PR~{\em (1):}
For every $s\in \mathbb{N}$ and for every $u \in Fact_\Lang$, 
we consider $\frac{p_\Lang(u,n_s)}{p_\Lang(n_s)}$, that is, the length of the
interval corresponding to $Cyl_\Lang(u)$ on the $x$-axis of $T_{n_s}$.
By a diagonal process on $\{\frac{p_{\Lang} (u,n_s)}{p_\Lang(n_s)}\}_{s \in
  \mathbb{N}}$, we obtain a sequence $Q=\{m_i\}_{i\in\mathbb{N}}$ such that, for
every $u$, $\lim_{i\rightarrow \infty} \frac{p_{\Lang} (u,m_i)}{p_\Lang(m_i)}$
exists, and we denote it by~$M_u$.
We then define $\mu$ on the cylinders of $\Lang$ by
$\mu(Cyl_\Lang(u))=M_u$. Note that $\mu(\Lang) = M_\epsilon = 1$.
The finite unions of the cylinders of $A^{\mathbb{N}}$ together with the empty
cylinder form a semiring of sets, thus the intersection of this semiring with
the cylinders of~$\Lang$ too, and $\mu$ is a premeasure on this semiring.
By Caratheodory's extension theorem, $\mu$ extends as a Borel probability
measure on~$\Lang$.

Note that $\lim_{i\rightarrow \infty} \frac{p_{\Lang}
  (u,m_i)}{p_\Lang(m_i)}=M_u$ implies that $\lim_{i\rightarrow \infty}
\frac{p_{\Lang} (u,m_i-1)}{p_\Lang(m_i-1)}=M_u$ too,   
that is, putting $d_i = \frac{p_{\Lang} (u,m_i)}{p_\Lang(m_i)} - \frac{p_{\Lang}
  (u,m_i-1)}{p_\Lang(m_i-1)}$, we have that $\lim_{i \rightarrow \infty} d_i=0$.
Indeed, first, since $\Lang$ is aperiodic minimal, $p_\Lang(m_i) >
p_\Lang(m_i-1)$ for every $i$,
thus
\vspace*{-0.2cm}
\[\begin{array}{llll}
d_i <   \frac{p_{\Lang} (u,m_i)}{p_\Lang(m_i-1)} -
            \frac{p_{\Lang} (u,m_i-1)}{p_\Lang(m_i-1)} 
      \le  |A| \frac{sp_{\Lang,l}(m_i-1)}{p_\Lang(m_i-1)},
\end{array}\]
\vspace*{-0.6cm}

\noindent
which, by Lemma~\ref{slopes_lem}, converges to~0. 
Next, since $p_{\Lang} (u,m_i) \ge p_{\Lang} (u,m_i-1)$, 
\vspace*{-0.1cm}
\[\begin{array}{llll}
d_i & \ge & \frac{p_{\Lang} (u,m_i-1)}{p_\Lang(m_i)} -
             \frac{p_{\Lang} (u,m_i-1)}{p_\Lang(m_i-1)} 
     =  p_{\Lang} (u,m_i-1) \frac{p_\Lang(m_i-1) - p_\Lang(m_i)}
                               {p_\Lang(m_i)p_\Lang(m_i-1)}  \\
    & \ge  & - p_{\Lang} (u,m_i-1) |A| 
             \frac{sp_{\Lang,l}(m_i-1)}{p_\Lang(m_i)p_\Lang(m_i-1)} 
      =   -|A| \frac{p_{\Lang} (u,m_i-1)}{p_\Lang(m_i)}
               \frac{sp_{\Lang,l}(m_i-1)}{p_\Lang(m_i-1)},
\end{array}\]
\vspace*{-0.4cm}

\noindent
which also converges to~0, again by Lemma~\ref{slopes_lem}. 

We now check the preservation of $\mu$ by $\sigma$.
Assume that the letter prolongations in $Fact_\Lang$ of $u$ to the left are
$a_1,..,a_h \in A$, where $h>0$ since $Fact_L$ is prolongable, $\Lang$ being
minimal.
Thus $\sigma^{-1}(Cyl_\Lang(u)) = \cup_{j=1}^h Cyl_\Lang(a_j u)$.
For all $m_i \in Q$, we also have
$\sum_{j=1}^h p_{\Lang}(a_ju, m_i) \geq p_{\Lang}(u, m_i-1)$,
and the difference between these two terms depends on the number of left special
factors in $Fact_{\Lang}(u, m_i-1)$ and on the number of their prolongations: 

\vspace*{-0.3cm}
{\small 
\EQ
0 \leq \sum_{j=1}^h p_{\Lang}(a_ju, m_i) - p_{\Lang}(u, m_i-1) \leq 
           |A|\; sp_{\Lang,l}(m_i-1).
\label{first_p_eq}
\EEQ
}
\vspace*{-0.4cm}

\noindent
Similarly, for each~$a_j$,  
$p_{\Lang}(a_ju, m_i) \geq p_{\Lang}(a_ju,m_i-1)$, 
and the difference between these terms depends on the number of right
special factors in $Fact_{\Lang}(a_ju,m_i-1)$: 

\vspace*{-0.3cm}
{\small
\EQ 
0 \leq \sum_{j=1}^h p_{\Lang}(a_ju, m_i) - \sum_{j=1}^h p_{\Lang}(a_ju, m_i-1)
\leq |A|\; sp_{\Lang,r}(m_i-1). 
\label{sec_p_eq}
\EEQ
}
\vspace*{-0.3cm}

\noindent
Reversing~(\ref{sec_p_eq}) and adding it to~(\ref{first_p_eq}) we get for all
$m_i \in Q$, 

\vspace*{-0.6cm}
{\small
$$ -|A|\; sp_{\Lang,r}(m_i-1) \le \sum_{j=1}^h p_{\Lang}(a_ju, m_i-1) -
p_{\Lang}(u, m_i-1) \le |A|\;sp_{\Lang,l}(m_i-1). 
$$ 
}
\vspace*{-0.5cm}

\noindent
We then divide the above terms by $p_\Lang(m_i-1)$, let $i$ go to infinity, and
using Lemma~\ref{slopes_lem} again, we obtain $\mu(Cyl_\Lang(u)) = \sum_{j=1}^h
\mu(Cyl_\Lang(a_j u))$.
Hence $\mu$ is an invariant measure for $\Lang$. 
Moreover, since $\Lang$ is aperiodic minimal, $\mu$ is nonatomic.

{\em (2):}
On the $x$-axis for each $u \in Fact_\Lang$, the length of the interval
corresponding to $Cyl_\Lang(u)$ converges to $M_u$ as $\lim_{i\rightarrow
  \infty} \frac{p_{\Lang} (u,m_i)}{p_\Lang(m_i)}$, and the same is true on the
$y$-axis for $\sigma(u)$ with $M_{\sigma(u)}$ as $\lim_{i\rightarrow \infty}
\frac{p_{\Lang} (\sigma(u),m_i-1)}{p_\Lang(m_i-1)}$.
Now, since $\Lang$ has zero topological entropy, and $\mu$ has all the requested
properties, Theorem~\ref{main1} applies to $\Lang$ measured by~$\mu$, so as to
obtain~$\Tiet$, with its associated partitions $PART_\Lang$ and
$PART_{\Lang,I}$.
The main components of $PART_\Lang$ are the cylinders $Cyl_\Lang(v^{(k)})$ (see
Lemma~\ref{fundlemma}),  
and following the construction of $T_{m_i}$, on the $x$-axis for each $k$, the
associated intervals $\Ik_{m_i}$ are the union of the intervals corresponding to
$Cyl_\Lang(v^{(k)}v)$, $v^{(k)}v \in Fact_\Lang(m_i)$, then converging to an
interval $\Ik$ of length $M_{v^{(k)}}$.
The same holds for the intervals $\Jk_{m_i}$ in the $y$-axis with
$Cyl_\Lang(\sigma(v^{(k)}v))$, respectively converging to intervals $\Jk$ of
length $M_{\sigma(v^{(k)})}$.
Note then that for every~$i$, $T_{m_i}$ is continuous, hence is an affine map on
each~$\Ik_{m_i}$.
Indeed, $\sigma$ is increasing on $Cyl_\Lang(\sigma(v^{(k)}))$, thus $T_{m_i}$
on $\Ik_{m_i}$ too. But if $T_{m_i}$ had a discontinuity, it would mean
that $\sigma^{-1}$ is not well-defined in this cylinder, contradicting the
definition of $PART_\Lang$.

Now, the left endpoint of $\Ik$ is $o_k=\sum_{v < v^{(k)}, v \in
  Fact_\Lang(|v^{(k)}|)} M_v$,
and the lower endpoint of~$\Jk$ is $o'_k= \sum_{v < \sigma(v^{(k)}), v \in
  Fact_\Lang(|v^{(k)}|-1)} M_{v}$.
By convergence, for every~$k$ and $\varepsilon > 0$, there is $j > 0$ such that
for all $i \ge j$, the endpoints of $\Ik_{m_i}$ are at distances $<\varepsilon$
respectively from $o_k$ and $o_k + M_{v^{(k)}}$.
Thus let~$T_{i, k, \varepsilon}$ be $T_{m_i}$ restricted to $\Ik_\epsilon=[o_k +
  \varepsilon, o_k + M_{v^{(k)}} - \varepsilon]$ for $i \ge j$,
on which it is affine with slope $\frac{p_{\Lang}(m_i)}{p_{\Lang}(m_i-1)}$.
These slopes converge to~$1$, and the endpoints of $\Jk_{m_i}$ converge
respectively to $o'_k$ and $o'_k + M_{\sigma(Cyl_\Lang(v^{(k)}))}$, hence $T_{i,
  k, \varepsilon}$ converges on $\Ik_\epsilon$.
Since $\varepsilon$ is arbitrary, $\{T_{m_i}\}_{i\in\mathbb{N}}$ converges on
all the $I^{(k)^o}$'s, i.e.~the $\Ik$'s interiors. We put $T$ as the
limit map, so that 
for each $x \in I^{(k)^o}$, 
$T(x) = o'_k + (x - o_k)$.
Now, in the proof of Lemma~\ref{veechlemma}, $T_\Lang(x) - x = K_{C^{(k)}}$,
where $C^{(k)}$ denoted $Cyl_\Lang(v^{(k)})$, and where
$K_{C^{(k)}} = \mu([w_{\Lang,min}, \sigma(w_{C^{(k)},min})]) -
                      \mu([w_{\Lang,min}, w_{C^{(k)},min}]$. 
By definition of $T$, the first term of this difference is $o'_k$ and the
second is~$o_k$, that is, 
$T_\Lang(x) - x = T(x) - x$ on $I^{(k)^o}$, 
thus by comparison with Lemma~\ref{veechlemma}, $I^{(k)^o} = I^{(k)^o}_\Lang$.
Since by Lemma~\ref{nullmeasure_infinspec_lem}, $\mu(\SP) = 0$,
the union of the~$I^{(k)^o}$'s has measure~$1$,
whence $T= T_\Lang$ almost everywhere.
\EPR

\noindent
Note that the map $T$ above can be extended to all of~$I$ by right-continuity,
so that $T$ can be made equal to $\Tiet$.

\subsection{Explicit Graph Examples}
\label{explicit_examples_sec}

\noindent
A dynamical system $(X,f)$ is said to be {\smbf uniquely ergodic} if there
exists only one $f$-invariant Borel probability measure on $X$.

\begin{Prop}\label{approx_ergod_uniq_prop} 
  Let $\Lang$ be an aperiodic minimal shift with linear complexity which is
  uniquely ergodic.  Then $\Tseq$ in Proposition~\ref{approx_prop}, defined as
  $\{T_n\}_{n>1}$, converges almost everywhere to $\Tiet$.
\end{Prop}
\PR Since $\Lang$ has linear complexity, its topological entropy is zero.  For
the same reason, there exists $r>0$ such that for all $n>1$, $p_\Lang(n) -
p_\Lang(n-1) \leq r$ \cite{Cas95}, hence $\lim_{n\rightarrow
  \infty}\frac{p_{\Lang}(n)}{p_{\Lang}(n-1)}=1$,
and Proposition~\ref{approx_prop} can be applied to~$\Lang$. 
Now, in the proof of this proposition, unique ergodicity means there can be only
one accumulation point for each $\{\frac{p_\Lang(u,n)}{p_\Lang(n)}\}_{n>1}$,
$u \in Fact_\Lang$. Thus these sequences converge without any subsequence
extraction, as does then $\Tseq$.
\EPR

\noindent
In order to illustrate the above proposition, the simplest non-trivial examples
are given by the primitive substitutive case (see p.~\pageref{primitive_def}):

\begin{Cor}\label{approx_ergod_uniq_subst_cor}
Let $\Lang$ be a shift associated with an aperiodic primitive substitution. Then
$\Tseq$ in Proposition~\ref{approx_prop}, defined as $\{T_n\}_{n>1}$, converges
almost everywhere~to~$\Tiet$.
\end{Cor} 
\PR Aperiodic primitive substitutions have associated shifts which are aperiodic
minimal with linear complexity~\cite{Pan84}, and uniquely ergodic~\cite{Que10}.
\EPR

\vspace*{-0.5cm}
\noindent
Note that in this primitive substitutive case, the limit measure of each
cylinder can also be computed from the substitution itself~\cite[Section
  5]{Que10}.

\begin{Example} \em
\label{thue-morse-ex-graph}
Reconsider the Thue-Morse substitution $\theta_{tm}$ (see
Examples~\ref{thue-morse-ex}, \ref{non-injectivity-ex},
and~\ref{infinite-discontinuities}).
Since $\theta_{tm}$ is aperiodic primitive,
Corollary~\ref{approx_ergod_uniq_subst_cor} applies to its associated shift
$\Lang_{tm}$.
In the following picture we show the graph of ${T}_{100}$ from the converging
sequence $\Tseq=\{T_n\}_{n>1}$ 
approximating the IET~$\Tiettm$ with $B_{\Lang_{tm}}$ as set of separation
points.
The intervals converging to their respective $I^{(k)}_{L_{tm}}$ of
$PART_{\Lang_{tm},I}$ are indicated for the first $k$'s, by the associated
$v^{(k)}$'s determining the corresponding cylinders of~$PART_{\Lang_{tm}}$ (see
Lemma~\ref{fundlemma} and Example~\ref{thue-morse-ex}).
As expected from the Examples~\ref{non-injectivity-ex}
and~\ref{infinite-discontinuities}, two accumulations of discontinuity points
can be observed, where non-injectivity holds.
Note first that, in this Thue-Morse case, the limit measure of each
involved cylinder can be directly obtained from~\cite{Dek92}.
Second, this example suggests existing connections with constructions of {\em
  infinite permutations} obtained from shifts associated with aperiodic
primitive substitutions~\cite{Mak09,AFP15f}.

\vspace*{0.1cm}
\begin{center} 
  \includegraphics[width=11.1cm]{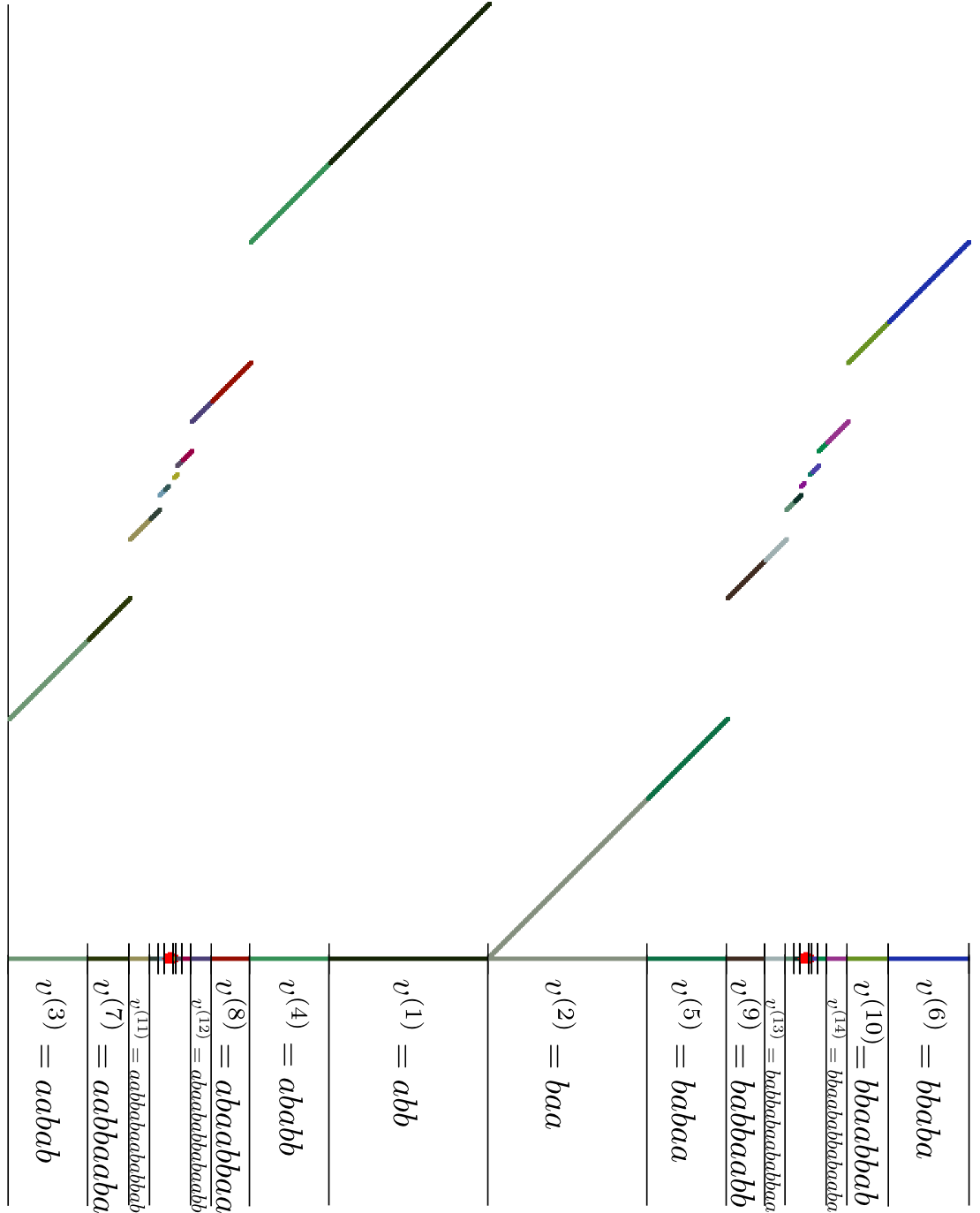}
\end{center}

\label{thue-morse_fig}
\end{Example}

\begin{Example}\em
\label{iet-graphs}
Here are other examples of approximated infinite IETs from linear complexity
shifts associated with aperiodic primitive substitutions. The pictures below are
the ${T}_{100}$ graphs from the convergent sequences
$\Tseq=\{T_n\}_{n>1}$ given by
Corollary~\ref{approx_ergod_uniq_subst_cor} for:
(i)~the {\em Tribonacci substitution}, i.e.~$\theta(a)=ab$,
$\theta(b)=ac$, $\theta(c)=a$; 
(ii)~the {\em Tetranacci substitution}, i.e.~$\theta(a)=ab$,
$\theta(b)=ac$, $\theta(c)=ad$, $\theta(d)=a$;
(iii)~the {\em Rudin-Shapiro substitution},
i.e.~$\theta(a)=ab$, $\theta(b)=ac$, $\theta(c)=db$, $\theta(d)=dc$; 
(iv)~the {\em Fibonacci substitution}, i.e.~$\theta(a)=ab$,
$\theta(b)=a$.
In this last example, $\Lang_\theta$ is known to be the shift given by the
natural coding of the minimal {\em finite} IET $T$ over two intervals, where
$B=\Disc=\{0,\frac{1}{\rho}\}$, and $\rho$ is the golden ratio~\cite{Fog02}.
Note then that in accordance with
Corollary~\ref{finite_iet_gives_the_same_finite_iet_cor}, the drawn graph is an
approximation of $T$'s graph and of $T_{\Lang_\theta}$'s graph too, while the
construction of $T_{\Lang_\theta}$ describes it with an infinite set
$B_{\Lang_\theta}$ of separation points.

\begin{center} 
\hspace*{-0.3cm}
\begin{tabular}{ccccccc}\label{iets_fig}
  \includegraphics[width=6.3cm]{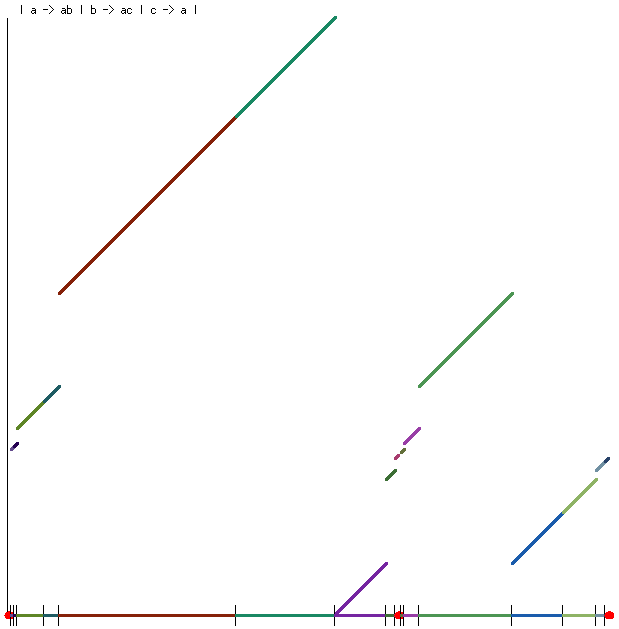} & \hspace*{0.6cm}
  \includegraphics[width=6.3cm]{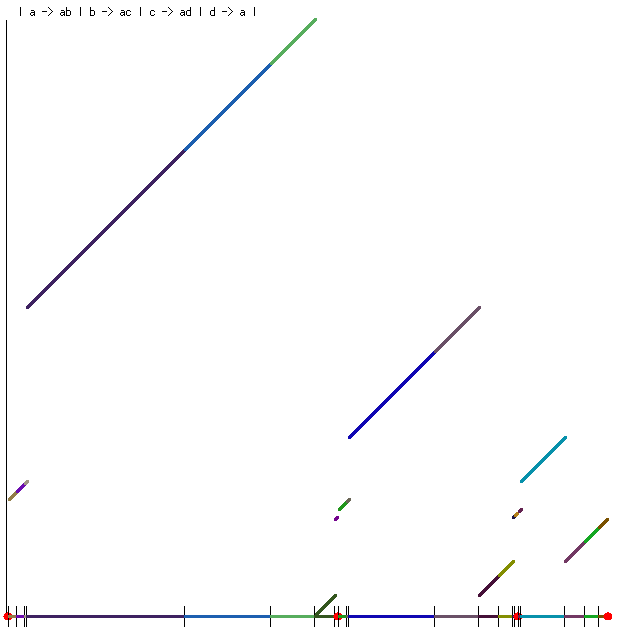} \\ 
  \em (i) Tribonacci & \hspace*{1cm} \em (ii) Tetranacci\\ \\ \\

  \includegraphics[width=6.3cm]{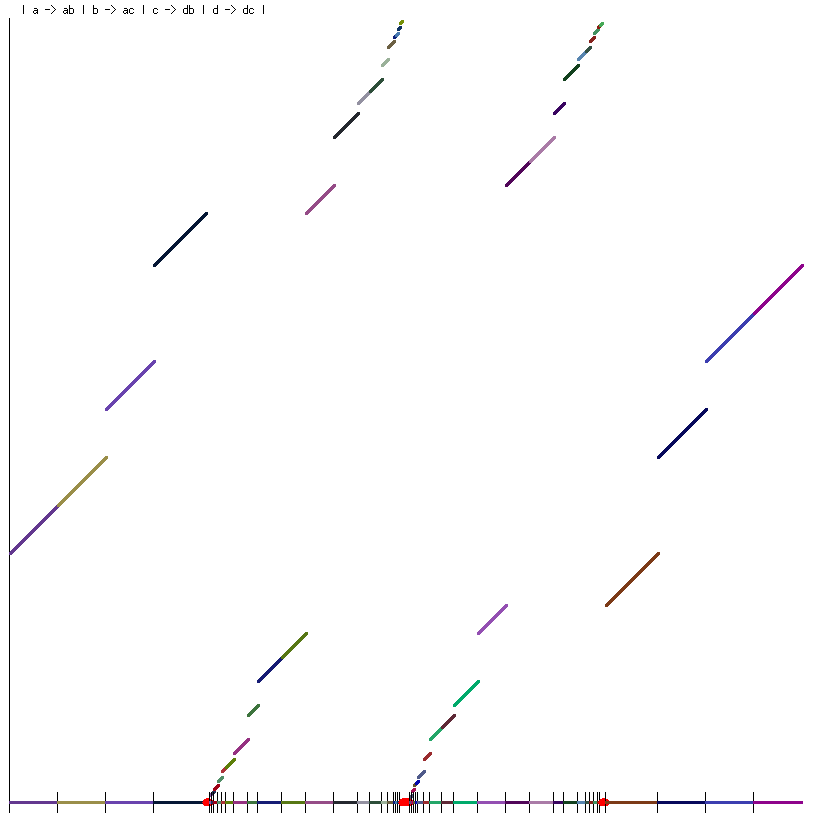} & \hspace*{0.6cm}
  \includegraphics[width=6.3cm]{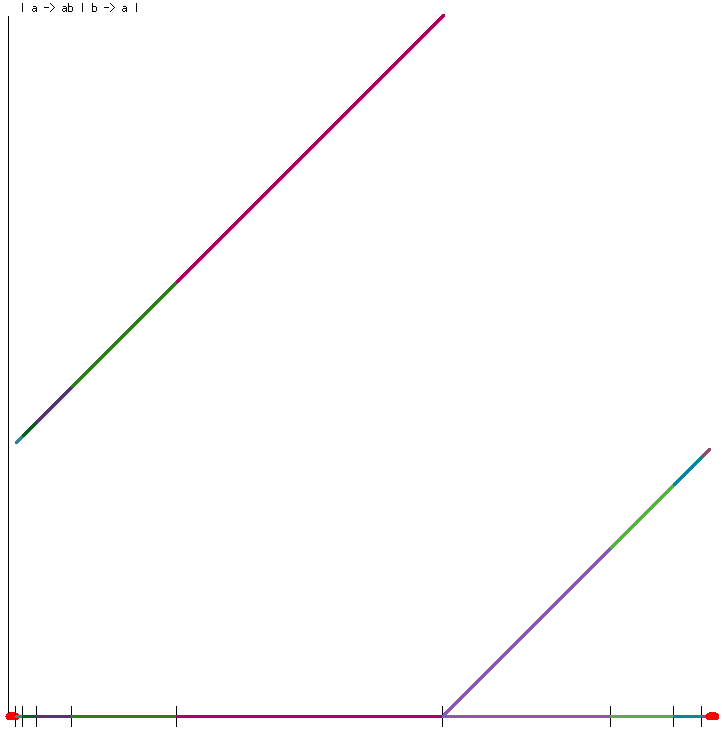} \\
  \em (iii) Rudin-Shapiro & \hspace*{1cm} \em (iv) Fibonacci
\end{tabular}
\end{center}
\end{Example}

\end{document}